\input amstex
\documentstyle{amsppt}
\magnification=\magstep1

\NoBlackBoxes
\TagsAsMath
\TagsOnRight

\pageheight{9.0truein}
\pagewidth{6.5truein}

\long\def\ignore#1\endignore{\par DIAGRAM\par}
\long\def\ignore#1\endignore{#1}

\ignore
\input xy \xyoption{matrix} \xyoption{arrow}
          \xyoption{curve}  \xyoption{frame}
\def\edge{\ar@{-}}
\def\dttdar{\ar@{.>}}

\def\dashedge{\ar@{--}}

\def\dshdar{\ar@{-->}}
\def\dropdown#1{\save+<0ex,-4ex> \drop{#1} \restore}

\def\loopNE{\ar@'{@+{[0,0]+(6,2)} @+{[0,0]+(10,10)}
@+{[0,0]+(2,6)}}}
\def\loopNW{\ar@'{@+{[0,0]+(-2,6)} @+{[0,0]+(-10,10)}
@+{[0,0]+(-6,2)}}}
\def\loopSW{\ar@'{@+{[0,0]+(-6,-2)} @+{[0,0]+(-10,-10)}
@+{[0,0]+(-2,-6)}}}
\def\loopSE{\ar@'{@+{[0,0]+(2,-6)} @+{[0,0]+(10,-10)}
@+{[0,0]+(6,-2)}}}

\def\loopNNE{\ar@'{@+{[0,0]+(4,2)} @+{[0,0]+(6,11)}
@+{[0,0]+(0,6)}}}
\def\loopSSW{\ar@'{@+{[0,0]+(-4,-2)} @+{[0,0]+(-6,-12)}
@+{[0,0]+(0,-6)}}}
\def\loopSSE{\ar@'{@+{[0,0]+(0,-6)} @+{[0,0]+(6,-11)}
@+{[0,0]+(4,-2)}}}
\endignore

\def\vsubseteq{\hbox{$\bigcup$\kern0.1em\raise0.05ex\hbox{$\tsize|$}}}
\def\smvsubseteq{\hbox{$\ssize\bigcup$\kern0.03em\raise0.05ex\hbox{$\ssize|$}}}

\def\la{{\Lambda}}
\def\lamod{\Lambda\text{-}\roman{mod}}

\def\AA{{\Bbb A}}

\def\PP{{\Bbb P}}
\def\SS{{\Bbb S}}

\def\NN{{\Bbb N}}

\def\hom{\operatorname{Hom}}
\def\aut{\operatorname{Aut}}

\def\top{\operatorname{top}}
\def\soc{\operatorname{soc}}
\def\topsocbd{\operatorname{\bold{Top-Soc(d)}}}

\def\orb{\operatorname{orbit}}

\def\dimvec{\operatorname{\underline{\dim}}}

\def\GL{\operatorname{GL}}

\def\Ext{\operatorname{Ext}}
\def\ext{\operatorname{ext}}

\def\autlap{\operatorname{Aut}_\la(P)}
\def\autlaphat{\operatorname{Aut}_{\lahat}(\Phat)}
\def\End{\operatorname{End}}

\def\udim{\operatorname{\underline{dim}}}

\def\Char{\operatorname{char}}

\def\C{{\Cal C}}
\def\D{{\Cal D}}

\def\F{{\Cal F}}

\def\M{{\Cal M}}

\def\U{{\frak U}}

\def\X{{\Cal X}}

\def\ba{\bold {a}}
\def\bb{\bold {b}}
\def\bc{\bold {c}}
\def\bC{\bold {C}}
\def\bd{\bold {d}}
\def\e{\bold {e}}
\def\bu{\bold {u}}
\def\bw{\bold {w}}

\def\Phat{\widehat{P}}
\def\Qhat{\widehat{Q}}

\def\bP{\bold{P}}
\def\ehat{\widehat{e}}

\def\chat{\widehat{c}}

\def\bdhat{\widehat{\bd}}
\def\buhat{\widehat{\bu}}
\def\lahat{\widehat{\la}}
\def\lahatmod{\lahat\text{-}\roman{mod}}

\def\alphahat{\widehat{\alpha}}

\def\Chat{\widehat{C}}

\def\ghat{\widehat{g}}

\def\Jhat{\widehat{J}}
\def\Mhat{\widehat{M}}

\def\Phat{\widehat{P}}
\def\Qhat{\widehat{Q}}

\def\Shat{\widehat{S}}
\def\That{\widehat{T}}

\def\xhat{\widehat{x}}

\def\bz{\bold{z}}

\def\modlad{\operatorname{\bold{Rep}}_{\bold{d}}(\Lambda)}

\def\Mod{\operatorname{\bold{Rep}}}

\def\rep{\operatorname{\bold{Rep}}}
\def\toptd{\operatorname{\bold{Rep}}^T_{\bold{d}}}

\def\laySS{\operatorname{\bold{Rep}} \SS}

\def\Hom{\operatorname{Hom}}

\def\grasstd{\operatorname{\frak{Grass}}^T_d}

\def\Gr{\operatorname{Gr}}
\def\grasstbd{\operatorname{\frak{Grass}^T_{\bold{d}}}}
\def\autlap{\aut_\Lambda(P)}

\def\unirad{\bigl(\aut_\la(P)\bigr)_u}

\def\grassSS{\operatorname{\frak{Grass}} \SS}

\def\grass{\operatorname{\frak{Grass}}}
 
 \def\biggrass{\operatorname{GRASS}}
 \def\biggrasslad{\operatorname{GRASS}_{\bd}(\la)}
 
 \def\biggrasstd{\operatorname{GRASS}^T_{\bd}}

\def\orbit{\operatorname{orbit}}

\def\term{\operatorname{end}}
\def\start{\operatorname{start}}
\def\Schu{\operatorname{Schu}}

\def\Schu{\operatorname{\Schu}}

\def\grasstd{\operatorname{\frak{Grass}}^T_d}

\def\grasstbd{\operatorname{\frak{Grass}}^T_{\bold d}}
\def\GRASS{\operatorname{GRASS}}
\def\biggrasstd{\GRASS^T_{\bold d}}
\def\biggrasstbd{\GRASS^T_{\bold d}}

\def\grassbd{\GRASS_{\bold d}(\Lambda)}

\def\Gr{\operatorname{Gr}}

\def\autlap{\aut_\Lambda(P)}
\def\autlat{\aut_\Lambda(T)}
\def\autlathat{\aut_{\lahat}(\That)}

\def\unirad{\bigl(\autlap \bigr)_u}

\def\grassSS{\operatorname{\frak{Grass}} \SS}
\def\biggrassSS{\GRASS\, \SS}

\def\grass{\operatorname{\frak{Grass}}}
\def\biggrass{\GRASS}

\def\orbit{\operatorname{orbit}}

\def\term{\operatorname{end}}
\def\start{\operatorname{start}}
\def\Schu{\operatorname{Schu}}

\def\AuReSm{{\bf 1}}
\def\BHT{{\bf 2}}
\def\BoHZtwo{{\bf 3}}
\def\CW{{\bf 4}}
\def\CBS{{\bf 5}}
\def\CKW{{\bf 6}}
\def\DoFl{{\bf 7}}
\def\EiSa{{\bf 8}}
\def\Ger{{\bf 9}}
\def\Gur{{\bf 10}}
\def\classifying{{\bf 11}}
\def\hier{{\bf 12}}
\def\moduli{{\bf 13}}
\def\GoHu{{\bf 14}}
\def\KacI{{\bf 15}}
\def\KacII{{\bf 16}}
\def\Kra{{\bf 17}}
\def\Mor{{\bf 18}}
\def\RiRuSm{{\bf 19}}
\def\Ros{{\bf 20}}
\def\Scho{{\bf 21}}
\def\Schro{{\bf 22}}

\topmatter

\title The geometry of finite dimensional algebras with vanishing radical square \endtitle

\rightheadtext{geometry for vanishing radical square}

\author Frauke M. Bleher, Ted Chinburg and Birge Huisgen-Zimmermann
\endauthor

\thanks While this research was carried out, the first author was partially supported by NSA grant H98230-11-1-0131, the second author by NSF grant DMS 110355, the third author by NSF grant DMS 0500961 and NSF grant DMS 0932078, while in residence at MSRI, Berkeley. \endthanks

\address Department of Mathematics, University of Iowa, Iowa City, IA 52242-1419 \endaddress 
\email frauke-bleher\@uiowa.edu \endemail

\address Department of Mathematics, University of Pennsylvania, Philadelphia, PA 19104-6395 \endaddress 
\email ted\@math.upenn.edu \endemail

\address Department of Mathematics, University of California, Santa
Barbara, CA 93106-3080 \endaddress
\email birge\@math.ucsb.edu \endemail

\abstract  Let $\la$ be a basic finite dimensional algebra over an algebraically closed field, with the property that the square of the Jacobson radical $J$ vanishes. We determine the irreducible components of the module variety $\Mod_\bd(\la)$ for any dimension vector $\bd$.  Our description leads to a count of the components in terms of the underlying Gabriel quiver.  A closed formula for the number of components when $\Lambda$ is local extends existing counts for the two-loop quiver to quivers with arbitrary finite sets of loops.  

For any algebra $\la$ with $J^2 = 0$, our criteria for identifying the components of $\Mod_\bd(\la)$ permit us to characterize the modules parametrized by the individual irreducible components.  Focusing on such a component, we explore generic properties of the corresponding modules by establishing a geometric bridge between the algebras with zero radical square on one hand and their stably equivalent hereditary counterparts on the other.  The bridge links certain closed subvarieties of  Grassmannians parametrizing the modules with fixed top over the two types of algebras.  By way of this connection, we transfer results of Kac and Schofield from the hereditary case to algebras of Loewy length $2$.  Finally, we use the transit of information to show that any algebra of Loewy length $2$ which enjoys the dense orbit property in the sense of Chindris, Kinser and Weyman has finite representation type.

\endabstract
   
\endtopmatter

\document

\head 1. Introduction \endhead

Let $\la$ be a basic finite dimensional algebra over an algebraically closed field $K$.  It is a fundamental problem to understand 
the finite dimensional representations of $\Lambda$.  One strategy for making headway in this direction is to determine the irreducible components of the varieties that parametrize the modules with fixed dimension vector, and to subsequently explore the generic structure of the modules corresponding to each of the components.  Investigations of this type were initiated by Kac in the hereditary case (where the parametrizing varieties  $\Mod_\bd(\la)$ are irreducible to begin with) and continued by Schofield, Crawley-Boevey, Schr\"oer, Carroll, Weyman, Babson, Thomas and the third author, among others.   The goal of the present paper is to advance this program for algebras with vanishing radical square.

Given any algebra $\la$ and a dimension vector $\bd = (d_1, \dots, d_n)$, two types of parametrizations of the $\la$-modules with this dimension vector are under consideration.  One is by the classical affine variety, $\modlad$, equipped with its $\GL(\bd)$-action. The other is by a projective variety, $\biggrasslad$, consisting of submodules of specified dimension of a suitable projective $\la$-module $\bP$, and endowed with the canonical $\aut_{\la}(\bP)$-action.  For a brief outline of these varieties and their relevant subvarieties, see Section 2 below. 

Without loss of generality, we may assume that $\la = KQ/I$ for a quiver $Q$ and an admissible ideal $I$ of the path algebra $KQ$.  In the hereditary case, that is for $I = 0$, the parametrizing varieties $\modlad$ and $\biggrasslad$ are always irreducible.  Work of Gabriel, Gelfand-Ponomarev, Kac, Ringel, Schofield, and many others led to a substantial theory relating the quiver $Q$ to generic features of the finite dimensional $KQ$-modules with fixed dimension vector.  Prominently, this work addresses decomposition properties shared by the modules in a suitable dense open subset of $\modlad$, as recorded by Kac's {\it canonical decomposition\/} of the dimension vector $\bd$.  Even when $\la$ has wild representation type, the Kac decompositions of dimension vectors are fully understood (see \cite{\KacI}, \cite{\KacII}, \cite{\Scho}).  

In attempts to extend results of this kind to the case where the ideal $I$ of relations is nonzero, one of the first hurdles encountered is the fact that the parametrizing varieties, $\modlad$ and $\biggrasslad$, split into a plethora of irreducible components in general.   Hence, one of the goals singled out as crucial early on  --   in the 1980's by Kraft at the latest (see \cite{\Kra})  --  is to determine the irreducible components of the parametrizing varieties in terms of the quiver $Q$ and the relations in  $I$.  This task is equivalent for the two types of varieties, affine and projective, but each offers methodological advantages over the other in certain situations; so it is advantageous to study them in parallel.  In the affine scenario, the task amounts to understanding the components of varieties consisting of finite sequences of matrices which satisfy certain relations.  Interest in such problems  is classical and precedes our present representation-theoretic focus; see e.g., \cite{\Ger}, \cite{\EiSa}, \cite{\Gur}.  

We next review some existing work regarding the nonhereditary case.  Results by Crawley-Boevey and Schr\"oer (see \cite{\CBS}) target the components $\C \subseteq \modlad$ whose representations are generically decomposable, relating them to components encountered for smaller dimension vectors. Here
``generically decomposable'' means that the modules corresponding to the points in some dense open subset of $\C$ are decomposable.  Moreover, in  \cite{\BHT}, Babson, Thomas and the third author have shown the component problem for the subvarieties $\biggrassSS\, $ of $\biggrasslad$ (cf\. Section 2) to be comparatively accessible; see Theorem 3.2.  Typically, this approach yields a finite collection of irreducible closed subsets of $\biggrasslad$ which includes the irreducible components of the big variety.  Hence the task of identifying the components of $\biggrasslad$ may be reduced to filtering them out of a larger collection of irreducible closed subsets.  Beyond these inroads into the problem, there are full descriptions of the irreducible components for a few special types of algebras, including the Gelfand Ponomarev algebras $K[X, Y]/ (XY, X^s, Y^t)$ (see \cite{\Schro}), the algebra $K[X,Y]/ (X^2, Y^2)$ (see \cite{\RiRuSm}), and some other special biserial algebras (see {\cite{\CW}).   

In this paper, we settle the component problem for the lowest interesting Loewy length, namely for vanishing radical square.  The main results in this connection are Theorem 3.6, Corollary 3.7, and Proposition 3.9.  As a consequence, we obtain a closed formula for the number of irreducible components of $\Mod_d(\la)$ in case $\la$ is local (Theorem 3.12); the number is given as a function of $d$ and the dimension $r$ of the radical.  Our formula subsumes existing computations for $r = 2$ (see \cite{\DoFl} and \cite{\Mor}).  Returning to arbitrary algebras with vanishing radical square, we proceed to study generic aspects of the representation theory supported by the individual components.  In particular, we explore the Kac decomposition of a dimension vector $\bd$ relative to an irreducible component of $\modlad$, as well as generic submodule lattices and tops.  The main results along this second line are Proposition 5.3 and Theorem 5.6.  They are first illustrated by way of local algebras, then applied to a problem of Chindris, Kinser, and Weyman.
\medskip

{\it More detailed outline of the paper. \/} Given a subset $\U$ of $\modlad$ or $\biggrasslad$, we refer to the modules that correspond to the points in $\U$ as the modules ``in" $\U$.  When $\U$ is an irreducible subvariety, the modules in $\U$ are said to {\it generically\/} have a property $(*)$ in case all modules in some dense open subset of $\U$ have property $(*)$.  Let $J$ denote the Jacobson radical of $\la$.  Two generically constant features of an irreducible component $\C$ of $\modlad$, which will play a pivotal role in the sequel, are the tops, $M/J M$,  and socles, $\soc M$, of the modules $M$ in $\C$.  We systematically identify isomorphic semisimple modules.   

In this outline {\it only\/}, we tacitly assume that $J^2 = 0$.

Section 2 assembles prerequisites concerning the parametrizing varieties. The main results of Section 3 characterize the irreducible components $\C$ of the varieties $\modlad$ (equivalently, those of the $\biggrasslad$).   These components are in one-to-one correspondence with the generic tops.  More precisely:  For any component $\C$ of $\modlad$, there exists a semisimple module $T$ with the property that $\C$ coincides with the closure of $\toptd$ in $\modlad$.  Here $\toptd$ is the subvariety of $\modlad$ that consists of the points corresponding to modules with top $T$; it is always irreducible, as is explained in Section 3.A.  In Corollary 3.7, we single out the generic tops, i.e., we determine those semisimple modules $T$ for which $\overline{\toptd}$ is maximal among the irreducible subvarieties of $\modlad$.  The theoretical description of the irreducible components of $\modlad$ has an algorithmic counterpart based on Proposition 3.9.  Indeed, the semisimple modules $T$ which identify the irreducible components of $\modlad$ can readily be determined from the quiver $Q$, as is illustrated by Example 3.11.  Theorem 3.12 applies the preceding results to the local case.

In Section 4, the $\la$-modules that belong to an irreducible component $\C$ with generic top $T$ are characterized.  This characterization, in turn, leads back to geometric information on the irreducible components of the affine, as well as the projective, parametrizing varieties.  In the affine case, the components of $\modlad$ are ``essentially" affine spaces, while the irreducible components of $\biggrasslad$ are ``close to" direct products of classical vector space Grassmannians.

In Section 5, we provide a geometric counterpart to the well-known fact that any algebra $\la$ with vanishing radical square is stably equivalent to a hereditary algebra $\lahat$. More specifically, we exhibit a strong link between certain parametrizing varieties over $\la$ and their analogues over $\lahat$ (Proposition 5.3). Subsequently we exploit this link to transfer information from the thoroughly studied hereditary case to the case $J^2 = 0$.  This connection appears to have been overlooked in the past, probably due to the fact that it cleanly surfaces only in the context of the varieties $\grasstbd$, projective counterparts of the varieties $\toptd$ (see Section 2).  We give two kinds of applications.  One is to Kac decompositions of dimension vectors and related generic features of the irreducible components (Theorem 5.6). It is illustrated in the local case (Section 6).  The other addresses a question which was raised and extensively studied by Chindris, Kinser and Weyman in \cite{\CKW} (Section 7).  They say that $\la$ has the {\it dense orbit property\/} if each irreducible component of any of the varieties $\modlad$ contains a dense $\GL(\bd)$-orbit.  Clearly, finite representation type implies the dense orbit property.  The converse is known to fail in general (see \cite{\CKW, Section 4}), but we will find it to be true whenever $J^2 = 0$ (Theorem 7.2). 
\smallskip

\noindent{\it Further conventions}:  The set $Q_0 = \{e_1, \dots, e_n\}$ of vertices of the quiver $Q$ will be identified with a full set of primitive idempotents of $\la$.   
The Cartesian product $(\NN_0)^n$ is equipped with the componentwise partial order 
$$(d_1,\dots,d_n) \le (d'_1,\dots,d'_n)  \iff  \ d_i \le d_i' \  \text{for}\ 1 \le i \le n.$$

For a $\la$-module $M$, we denote by $\top M$ its top, and by $\SS(M)$ its radical layering $(J^lM/J^{l+1}M)_{0 \le l \le L}$, where $L$ is maximal with $J^L \ne 0$.  We continue to identify isomorphic semisimple modules, unless we wish to distinguish among different embeddings.  The radical layerings are examples of {\it semisimple sequences\/}, i.e., sequences $\SS = (\SS_0, \dots, \SS_L)$ of semisimples $\SS_l$ in $\lamod$.  The {\it top of $\SS$\/} is $\SS_0$, the {\it dimension vector of $\SS$\/}, denoted $\udim\, \SS$, is the dimension vector of the direct sum $\bigoplus_{0 \le l \le L} \SS_l$ of the entries.  Moreover, a {\it top element\/} of $M$ is any element $z \in M \setminus JM$ which is normed by one of the primitive idempotents in $Q_0$, i.e., $e_i z = z$ for some $i\le n$.  A {\it full sequence of top elements of $M$\/} is any generating set of $M$ consisting of top elements which are $K$-linearly independent modulo $JM$.

\head 2. Prerequisites on parametrizing varieties \endhead

We start with an overview, Diagram 2.1, of the relevant varieties (slightly modifying the notation of \cite{\hier}, where more detail can be found). Then we recall how information is transferred among these varieties along the horizontal two-way arrows.  Let $\bd = (d_1, \dots, d_n)$ be a dimension vector with $d = |\bd| = \sum_i d_i$ and $\bP = \bigoplus_{1 \le r \le d} \la \bz_r$ a projective $\la$-module such that $\udim (\top \bP) = \bd$; here $\bz_1, \dots, \bz_d$ is a full sequence of top elements of $\bP$.  In other words, $\bP$ is a projective cover of $\bigoplus_{1 \le i \le n} S_i^{d_i}$.  For a semisimple $\la$-module $T$ with $\udim T = (t_1, \dots, t_n)$ and total dimension $t = \sum_i t_i$, we fix a projective cover $P = \bigoplus_{1 \le r \le t} \la z_r$ of $T$; here the $z_r$ constitute a full sequence of top elements of $P$.  Observe that any $\la$-module with dimension vector $\bd$ arises as a factor module of $\bP$ and every module with top $T$ arises as a factor $P/C$ for a suitable submodule $C \subseteq JP$.   The motivation for introducing the ``small scenario", based on $P$ instead of $\bP$, lies in the fact that all problems concerning the generic behavior of modules can be resolved in the small, much more manageable, setting (see Section 5 for illustration).    

\goodbreak \midinsert
$$\xymatrixrowsep{0.2pc}\xymatrixcolsep{3pc}
\xymatrix{
\ar@{}[r]|{\txt{``big'' scenario}}   &&\txt{``small'' scenario}  \\
\boxed{\grassbd} \dropdown{\txt{(projective)}} \ar@{<->}[r] &\boxed{\modlad}
\dropdown{\txt{(affine)}} \\  \\
\smvsubseteq &\smvsubseteq \\
\boxed{\biggrasstbd} \dropdown{\txt{(quasi-projective)}} \ar@{<->}[r]
&\boxed{\toptd} \dropdown{\txt{(quasi-affine)}}
\ar@{<->}[r] &\boxed{\grasstbd} \dropdown{\txt{(projective)}} \\  \\
\smvsubseteq & \smvsubseteq & \smvsubseteq \\
\boxed{\biggrassSS} \dropdown{\txt{(quasi-projective)}} \ar@{<->}[r]
&\boxed{\laySS} \dropdown{\txt{(quasi-affine)}}
\ar@{<->}[r] &\boxed{\grassSS} \dropdown{\txt{(quasi-projective)}} \\  \\
}$$

\medskip
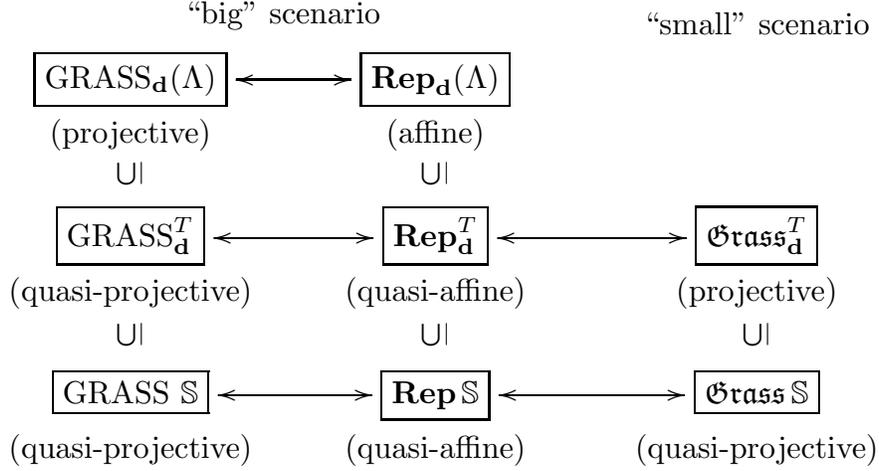
\centerline{Diagram 2.1}
\endinsert  

Here the horizontal double arrows point to the strong geometric correspondences spelled out in Proposition 2.1.  We follow with definitions of the entries.
The varieties in the left-most column of Diagram 2.1 are subvarieties of the classical vector space Grass\-mann\-ian $\Gr(\dim \bP - d, \bP)$ consisting of the $(\dim\bP - d)$-dimensional $K$-subspaces of $\bP$;  those in the right-most column are subvarieties of the Grassmannian $\Gr(\dim P - d, P)$.  We introduce the  displayed varieties from top to bottom, moving from right to left in each row. Let  $T$  be a semisimple $\la$-module with dimension vector $\le\bd$, and $(\SS_0, \dots, \SS_L)$ a semisimple sequence with top $T$ and dimension vector $\bd$. 

\roster
\item $\modlad$: This is the classical affine variety parametrizing the left $\la$-modules with dimension vector $\bd$, namely

\noindent$ \bigl\{ (x_\alpha)_{\alpha \in Q_1} \in \prod_{\alpha \in Q_1} \Hom_K \bigl(K^{d_{\text{start}(\alpha)}},\, K^{d_{\text{end}(\alpha)}}\bigr) \mid \text{the\ } x_{\alpha} \ \text{satisfy all relations in}\ I \bigr\},$

\noindent where $Q_1$ is the set of arrows of the quiver $Q$.
 As usual, we endow $\modlad$ with the conjugation action of $\GL(\bd) := \GL_{d_1}(K) \times \cdots\times \GL_{d_n}(K)$, the orbits of which are in bijective correspondence with the isomorphism classes of modules that have dimension vector $\bd$.

\item $\grassbd$: Set $a := \dim\bP-d$. Then $\grassbd$ is the closed subvariety of the Grassmannian $\Gr(a,\bP) \subseteq \PP(\la^a\bP)$ consisting of those $a$-dimensional $K$-subspaces $C \subseteq \bP$ which are $\la$-submodules of $\bP$ such that $\udim(\bP/C)=\bd$. In particular, $\grassbd$ is a projective variety. Note, moreover:  The linear algebraic group $\aut_{\la}(\bP)$ acts mor\-phically on  $\grassbd$, and the orbits of this action are, in turn, in 1--1 correspondence with the isomorphism classes of left $\la$-modules having dimension vector $\bd$.

\item $\grasstbd$ denotes the closed subvariety of the Grassmannian $\Gr(\dim P - d, P)$ consisting of those $K$-subspaces $C\subseteq P$ which are $\la$-submodules of $JP$, with the additional properties that $P/C$ has top $T$ and $\udim(P/C)= \bd$.  Clearly, $\grasstbd$ carries a morphic action of the smaller automorphism group $\autlap$.  This time, the orbits of our group action are in 1--1 correspondence with the isomorphism classes of left $\la$-modules with top $T$ and dimension vector $\bd$.  Without loss of generality, we may identify $\bz_r$ with $z_r$ for $r \le t$, which yields an  embedding $P \subseteq \bP$ via $C \mapsto C \oplus \bigoplus_{t +1 \le r \le d} \la \bz_r$.  This embedding makes $\grasstbd$ a closed subvariety of $\grassbd$ which, however, fails to be closed under the action of $\aut_\la(\bP)$ in general. 

\item $\toptd$: This is the locally closed subvariety of $\modlad$ consisting of the points corresponding to modules with top $T$. It is clearly closed under the $\GL(\bd)$-action.

\item $\biggrasstbd$ is the locally closed subvariety of $\grassbd$ consisting of all points  that correspond  to modules with top $T$. Clearly,  the embedding $\grasstd \hookrightarrow \grassbd$ under (3) places $\grasstd$ into $\biggrasstbd$.  We point out that $\biggrasstbd$ is stable under the $\aut_{\la}(\bP)$-action; in fact $\biggrasstbd$ is the closure of $\grasstd$ under the action of the big automorphism group.  

\item $\grassSS$, $\laySS$, $\biggrassSS$: Each of these is the locally closed subvariety of the variety shown above it that consists of the points parametrizing the modules with radical layering $\SS$. Evidently, in each case, the mentioned subvariety is closed under the pertinent group action, that of $\autlap$, $\GL(\bd)$, and $\aut_{\la}(\bP)$, respectively.  Observe that $\biggrassSS$ coincides with the closure of  $\grassSS$ under  the $\aut_{\la}(\bP)$-action.  The varieties $\grassSS$ and $\biggrassSS$ are very similar geometrically in that they have open affine covers, the patches of which differ only by a (large) direct factor $\AA^m$. 
\endroster

The various settings were connected by Bongartz and the third author in \cite{\BoHZtwo, Proposition C}.  We only quote the transfer of information that is relevant here.

\proclaim{Proposition 2.1. Transfer between the affine and projective settings}   
\smallskip

{\bf (I) $\grassbd$ versus $\modlad$:}  Consider the one-to-one correspondence between the orbits of $\grassbd$ and $\modlad $ assigning to any orbit $\aut_\la(\bP).C \subseteq \grassbd$ the orbit
$\GL(\bd).x \subseteq \modlad $ that represents the same isomorphism class of
$\la$-modules.  This correspondence extends to an inclusion-pre\-serv\-ing bijection
$$\Phi: \{ \aut_\la(\bP)\text{-stable subsets of\ } \grassbd \} \rightarrow
\{\GL(\bd)\text{-stable subsets of\ } \modlad \}$$  
which preserves and reflects
openness, closures and irreducibility.  In particular, $\Phi$ takes $\biggrasstbd$ to $\toptd$ and $\biggrassSS$ to $\laySS$. 
\smallskip

{\bf (II) $\grasstd$ versus $\toptd$:}  The one-to-one correspondence which
pairs any orbit $\autlap.C$ of $\grasstd$ with the orbit
$\GL(\bd).x$ of $\toptd$  representing the same isomorphism class of $\la$-modules extends to an inclusion-preserving bijection
$$\phi: \{ \autlap\text{-stable subsets of\ } \grasstd \} \rightarrow
\{\GL(\bd)\text{-stable subsets of\ } \toptd \}$$  
which, once again, preserves and reflects
openness, closures and irreducibility.

\endproclaim

\proclaim{Observation 2.2} {\rm{(See \cite{\BHT, Section 2}) }}  Each irreducible component of any of the varieties $X$ in rows one and two of Diagram 2.1 is among the closures in $X$ of  the irreducible components of the subvariety below it.  However, the set of these closures typically also contains irreducible subvarieties of $X$ which fail to be components of $X$.  \endproclaim 

It will turn out that the interplay among tops and socles of the modules in the irreducible subvarieties of $\modlad$ is key to identifying the components in case $J^2 = 0$.  The following observation concerning generic tops and socles is not new.  It rests on the well-known fact that, for any $N \in \lamod$ and dimension vector $\bd$, the functions $\modlad \rightarrow \NN_0$,
$$x \mapsto  \dim \Hom_\la(N, M(x)) \ \ \ \text{and} \ \ \ x \mapsto  \dim \Hom_\la(M(x), N)$$
are upper semicontinuous.   Here $M(x)$ denotes the $K$-space $\bigoplus_{1 \le i \le n} K^{d_i}$ with the $\la$-module structure induced by the linear maps $x_{\alpha}: K^{d_i} \rightarrow K^{d_j}$, where $\alpha$ is an arrow from $e_i$ to $e_j$.  

\proclaim{Observation 2.3}  Suppose that $\U$ is an irreducible subvariety of $\modlad$.  Then there exists a dense open subset $\U' \subseteq \U$ such that, for all $x \in \U'$, the dimension vector of the top of $M(x)$ is the unique minimal element of $\{ \udim \top M(y) \mid y \in \U\}$, and the dimension vector of the socle of $M(x)$ is the unique minimal element of $\{ \udim \soc M(y) \mid y \in \U\}$.

In short, the generic dimension vectors of tops and socles of the modules in $\U$ are the minimal ones attained on $\U$.  \qed
\endproclaim

\head 3.  Identifying the irreducible components of the varieties $\modlad$ and $\biggrass_{\bd}(\la)$ for $J^2 = 0$. \endhead

An initial source of irreducible components of $\modlad$ arises as an immediate consequence of Observation 2.3.  The following twin pair of observations applies to {\it any\/} finite dimensional algebra:
 
\proclaim{Observation 3.1} Let $\bd$ be a dimension vector, and let $T$, $U$ $\in \lamod$ be semisimple.  By $\Mod^{\bd}_{U}$ we denote the locally closed subvariety of $\modlad$ consisting of the points representing modules with  socle $U$. 
\smallskip

 {\bf (a)} If $\udim T$ is minimal among the dimension vectors $\udim (\top M)$, where $M$ traces the modules with dimension vector $\bd$, then the closure in $\modlad$ of any irreducible component of $\toptd$ is an irreducible component of $\modlad$.
 \smallskip
 
 {\bf (b)}  If $\udim U$ is minimal among the dimension vectors $\udim (\soc M)$, where $M$ traces the modules with dimension vector $\bd$, then the closure in $\modlad$ of any irreducible component of $\Mod^{\bd}_{U}$ is an irreducible component of $\modlad$. \qed
 \endproclaim
 
 However, even in case $J^2 = 0$, one usually misses the majority of irreducible components of $\modlad$ by restricting to those provided by Observation 3.1.
To make further headway in this case, we incrementally tighten our restrictions on the algebra $\la$.
\medskip

\noindent {\bf{(A)  The component problem for truncated path algebras}} 
\smallskip

Every finite dimensional basic algebra with $J^2 = 0$, as well as any hereditary algebra, is a {\it truncated path algebra\/} in the following sense:  There exist a quiver $Q$ and an integer $L \ge 1$ such that $\la = KQ/ I$, where $I$ is the ideal generated by all paths of length $L+1$.   In this situation, the key ingredients in the quest for the irreducible components of the $\modlad$ are the semisimple sequences $\SS$ with dimension vector $\bd$.

\proclaim{Theorem 3.2} \cite{\BHT, Theorem 5.3}  Let $\la$ be a truncated path algebra.

Each of the nonempty subsets $\laySS$ of $\modlad$ is an irreducible, smooth,  unirational subvariety of $\modlad$.  Moreover, every irreducible component of $\modlad$ arises as the closure $\overline{\laySS}$ of some subvariety $\laySS$.  
\endproclaim 

Suppose $\SS$ is an arbitrary semisimple sequence with dimension vector $\bd$ over a truncated path algebra $\la$ such that $\laySS \ne \varnothing$.  Clearly, the irreducible subvariety $\overline{\laySS}$ may fail to be maximal among the irreducible subsets of $\modlad$.  Indeed, if we take $\SS$ to be of the form $(\SS_0, 0, \dots, 0)$ with $\udim \SS_0 = \bd$, then $\overline{\laySS} = \laySS$ is contained in every orbit closure of $\modlad$. 

Theorem 3.2 pinpoints a special feature of truncated path algebras.  Even for a monomial algebra $\la$, there are typically multiple irreducible components of $\modlad$ whose modules have the same generic radical layering:  For a small example, consider the quiver $Q$ with $3$ vertices, and $4$ arrows, $\alpha_1, \alpha_2$ from $e_1$ to $e_2$ and $\beta_1, \beta_2$ from $e_2$ to $e_3$.  Let $I$ be the ideal generated by the single relation $\beta_2\alpha_2$.  For $\bd = (1, 1, 1)$, the variety $\modlad$ has two irreducible components, both of which have generic radical layering $\SS = (S_1, S_2, S_3)$; generically, the modules in these components have graphs of the form
$$\xymatrixrowsep{1.4pc}\xymatrixcolsep{4pc}
\xymatrix{
1 \edge@/_/[d]_{\alpha_1} \edge@/^/[d]^{\alpha_2} &&1 \edge[d]^{\alpha_1}  \\
2 \edge[d]_{\beta_1} &\txt{and} &2 \edge@/_/[d]_{\beta_1} \edge@/^/[d]^{\beta_2}  \\
3 &&3
}$$

 In light of Theorem 3.2, the task of determining the irreducible components of $\modlad$ for a truncated path algebra $\la$ translates into the problem of sifting out those  semisimple sequences $\SS$ which have the property that the closure $\overline{\laySS}$ in $\modlad$ is {\it maximal\/} irreducible.  This is a very challenging problem in general.  However, it becomes quite accessible if one restricts to algebras with Loewy length $2$.  
 
\medskip

\noindent {\bf{(B)  Focus on algebras with vanishing radical square}}
\smallskip

We adopt the following blanket hypothesis for the remainder of this section:  $J^2 = 0$.
 \smallskip

In this case, the relevant semisimple sequences $\SS$ have only two entries, $\SS = (\SS_0, \SS_1)$, and for fixed dimension vector $\bd$, we obtain  
$$\laySS = \toptd, \ \ \ \  \text{where} \ \ T = \SS_0.$$ 
Consequently, determining the radical layerings $\SS$ that are generic for the irreducible components of $\modlad$ amounts to  pinning down the tops $T$ that are generic for the components.  The following terminology will be convenient in this connection.

\definition{Definition  and Comments 3.3} Let $\bd$ be a dimension vector, and let $T \in \lamod$ be semisimple.
\smallskip

$\bullet$ $T$ will be called {\it realizable with respect to $\bd$\/} in case $\toptd \ne \varnothing$, i.e., in case $\bd - \udim T \le \udim JP$, where $P$ is the projective cover of $T$.
\smallskip

$\bullet$  $T$ will be called a {\it generic top of $\modlad$\/} if $\overline{\toptd}$ is an irreducible component of $\modlad$.  In that case, the generic socle of  the modules in $\overline{\toptd}$ will also be referred to as a {\it generic socle of $\modlad$\/}. 
\smallskip

$\bullet$ The partially ordered set $\topsocbd$:  Define the following set of pairs
$$\topsocbd = \{ (\udim \top M,\, \udim \soc M) \mid M\in \lamod\, \ \text{and}\ \, \udim M = \bd \}.$$
It is equipped with the componentwise partial order on pairs in $(\NN_0)^n \times (\NN_0)^n$, based on our partial order of $(\NN_0)^n$ (cf\. conventions).
\enddefinition

The upcoming remarks are straightforward. 

\proclaim{Observation 3.4}  Let $\SS = (\SS_0, \SS_1)$ be a semisimple sequence with dimension vector $\bd$.
\smallskip

{\bf (a)}  If $M \in \lamod$ has radical layering $\SS$, then $\SS_1 \subseteq \soc M$, and equality holds precisely when $M$ has no simple direct summands.  In fact, suppose $M$ is decomposed in the form $M = M_0 \oplus M_1$, where $M_0$ is without simple direct summand and $M_1$ is semisimple; then $\soc  M = JM_0 \oplus M_1$ and  $J M = JM_0 = \soc M_0$. 

Thus, $\udim \top M + \udim \soc M = \bd$ if and only if $M$ has no simple direct summand.
\smallskip

{\bf (b)} The finite set $\topsocbd$ can be algorithmically determined from $Q$ and $\bd$.  Instead of spelling out a formal procedure, we exemplify it below.  
\smallskip

{\bf (c)} For every pair $(\ba\, , \bb) \in \topsocbd$, the sum $\ba + \bb$ is bounded from below by $\bd$.  In particular, $(\ba\, , \bb)$ is minimal in $\topsocbd$ whenever $\ba + \bb = \bd$.  We will find the converse to be true in many interesting instances {\rm {(}}see {\rm Theorem 3.12)}.  But in general the sums $\ba + \bb$, as $(\ba, \bb)$ traces the minimal elements of $\topsocbd$, will vary; see the example preceding {\rm Proposition 3.9}, or {\rm Example 3.11}. \qed
\endproclaim

For later reference, we compute the set $\topsocbd$ in a specific instance.

\definition{Example 3.5}  Let $\la = KQ/I$, where $Q$ is the quiver with a single vertex and $3$ loops, and $I = \langle \text{all paths of length}\ 2 \rangle$.  Take $\bd = d = 13$.  Moreover, denote by $S$ the unique simple left $\la$-module.  Given that the smallest possible top of a module with dimension $d$ is $S^4$ and the largest is $S^{13}$, we find  
$\topsocbd$ to consist of the following pairs:  $(4, 9)$, $(4, 10)$, $(5, 8)$, $(5, 9)$, $(5, 10)$, $\allowmathbreak$ $(6, 7), (6,8), \dots, (6, 10)$, $(7, 6), (7, 7), \dots, (7, 11)$, $(8, 5), (8, 6), \dots, (8, 11)$, $(9, 4), (9,5), \dots,$ $\allowmathbreak$  $(9, 11)$, $(10, 4), (10, 5), \dots, (10, 12)$, $(11, 7),(11, 8), \allowmathbreak  \dots, (11, 12)$, $(12, 9), (12, 10), \dots, (12, 12)$, $(13, 13)$.  In particular, the minimal pairs are precisely those pairs $(a, b)$ for which $4 \le a \le 9$ with $a + b = 13$.

We present sample arguments indicating the method:  Suppose that $\dim M = 13$, and let $M = M_0 \oplus M_1$ be as in Observation 3.4(a) .   Clearly, $\dim \top M \ge 4$, since for $\dim \top M \le 3$, we obtain $\dim JM \le 9$, which would entail $\dim M < 13$. Now we focus on two specific values of  $\dim \top M$.  First consider the case $\dim \top M = 7$.  Then $\dim JM = \dim JM_0 = 6$, and the possible choices for $\dim  \top M_0$ are precisely $S^7, S^6, \dots, S^2$, with $S^2$ the smallest due to $JM_0 = \soc M_0 =  S^6$.  Since $\top M =  M_0 \oplus \top M_1$ and $\soc M =  \soc M_0 \oplus M_1$, the corresponding socles of $M$ are $0 \oplus S^6, S \oplus S^6,  \dots, S^5 \oplus S^6$.  We demonstrate the dual argumentation for $\dim \top M = 11$.   In that case, the injective envelope of $JM = S^2$ has top dimension $6$ whence $\dim \top M_0 \le 6$.  So $M_1$ contains $S^5$, and the smallest occurring socle is $S^7$;  on the other extreme, the smallest possible top of $M_0$ is $S$, in which case $M_1 = S^{10}$ and $\soc M = S^{12}$. 
\enddefinition

The main results of this subsection give two equivalent characterizations of the {\it generic tops\/} of $\modlad$.  The second, obtained by way of the first, facilitates the assembly of the list of generic tops from the quiver $Q$. These characterizations, in turn, lead to an explicit description of the irreducible components of $\modlad$.   

\proclaim{Theorem 3.6}  Suppose that $J^2 = 0$.  Let $T_1$, $T_2  \in \lamod$ be semisimple modules, both realizable with respect to $\bd$.  Moreover, let $U_i$  be the generic socle of the modules in $\rep^{T_i}_{\bd}$ for $i = 1,2$; in particular, $(\udim T_i, \udim U_i) \in \topsocbd$.  Then
$$\rep^{T_1}_{\bd} \subseteq \overline{\rep^{T_2}_{\bd}} \ \ \iff \ \ (\udim T_1, \udim U_1) \ge (\udim T_2, \udim U_2).$$
\endproclaim 

\demo{Proof}    
The implication ``$\implies$" is clear from Observation 2.3.  For the converse, suppose 
$$(\udim T_1, \udim U_1) > (\udim T_2, \udim U_2)\ \ \ (\dagger).$$
Then $T_1 \supsetneqq T_2$ by the definition of the $U_i$, say $T_1 = T_2 \oplus V$, where $V$ is a nonzero semisimple module.  Let $X$ be any module in $\rep^{T_1}_{\bd}$ with $\soc X = U_1$.    Clearly, it suffices to show that $X$ belongs to $\overline{\rep^{T_2}_{\bd}}$ under these circumstances.  

Realizability of $T_2$ with respect to $\bd$ and the choice of $U_2$ permit us to pick a  module $Y$ in $\rep^{T_2}_{\bd}$ with $\soc Y = U_2$.   From $(\dagger)$ we deduce that $\soc X = U_1 \supseteq U_2 = \soc Y$.  Moreover, in light of the equalitites $\udim T_2 + \udim V + \udim JX = \udim T_1 + \udim JX = \udim X  =  \udim Y = \udim T_2 + \udim JY$, we infer
$$JX \oplus V =  JY  \subseteq \soc Y \subseteq \soc X$$
by our convention of identifying isomorphic semisimple modules.  This means that the annihilator of $J$ in $X$ contains $JX \oplus V$.  Therefore $X  \cong X' \oplus V$, where $\top X' = T_2 = \top Y$ and $JX = JX'$.   Let $P$ be a projective cover of $X'$, and thus also of $Y$.  It is harmless to assume $X' = P/C'$ for a suitable submodule $C' \subseteq JP$.  
Write $\udim V = (v_1, \dots, v_n)$, $\underline{\dim}\,  C' = (w_1, \dots, w_n)$ and, for each $i \in \{1, \dots, n\}$, pick a  $K$-basis $c_{i1}, \dots, c_{i, w_i}$ for $e_i C'$.  Due to semisimplicity of $JP$, each of the one-dimensional spaces $K c_{ij}$ is a direct summand of $JP$ which is isomorphic to $S_i$.  In light of $J Y = JX' \oplus \bigoplus_{1 \le i \le n} S_i^{v_i} \subseteq JP$, we infer that $v_i \le \dim e_i C'  = w_i$.  This allows us to define a submodule $C$ of $C'$ as follows: 
$$C = \bigoplus_{1 \le i \le n}\ \bigoplus_{v_i + 1 \le  j \le w_i} K c_{ij}.$$ 
Note that $\udim C = (w_1 - v_1, \dots, w_n - v_n)$.  Setting $Z = P/ C$, we thus find: $\udim Z = \bd$, $\top Z = T_2$, and $C'/C$ is a semisimple submodule of $Z$ with dimension vector $(v_1, \dots, v_n)$.  The short exact sequence
$$0 \longrightarrow \ C' / C\ \longrightarrow \ Z = P/ C \longrightarrow \  X' = P /C'  \longrightarrow \ 0$$
now shows that $Z$ degenerates to $X' \oplus \bigl(C'/C \bigr) = X' \oplus  \bigoplus_{1 \le i \le n} S_i^{v_i} = X$.  This places $X$ into the closure $\overline{\rep^{T_2}_{\bd}}$ and completes the argument.  \qed 
\enddemo

\proclaim{Corollary 3.7} Suppose $J^2 = 0$, and let $T \in \lamod$ be semisimple.  Then the closure of $\toptd$ in $\modlad$ is an irreducible component of $\modlad$ if and only if there exists a module $M$ with top $T$ and dimension vector $\bd$ such that the pair $(\dimvec \top M,\, \dimvec \soc M)$ is minimal in the set $\topsocbd$.
\smallskip

If $M_1, \dots, M_m$ are modules with dimension vector $\bd$ such that the pairs
$$(\dimvec \top M_1,\, \dimvec  \soc M_1), \dots, (\dimvec \top M_m,\, \dimvec \soc M_m)$$
are the distinct minimal elements of $\topsocbd$,
then $\modlad$ has precisely $m$ distinct irreducible components, namely the closures of the subvarieties $\toptd$, where $T$ traces $\top M_1, \dots, \top M_m$.  In particular, $\top M_1, \dots, \top M_m$ are the generic tops of $\modlad$, and $\soc M_1, \dots, \soc M_m$ are the generic socles of $\modlad$.
\endproclaim 

\demo{Proof}  By Theorem 3.2 and the opening comments of Subsection B, each irreducible component of $\modlad$ is of the form $\overline{\toptd}$ for some semisimple module $T$ which is realizable with respect to $\bd$.  Hence the first claim follows from Theorem 3.6. 
To justify the remaining assertions, set $T_j = \top M_j$ and $\C_j = \overline{\rep^{T_j}_{\bd}}$.  Clearly $T_j$ is then the generic top of the modules in $\C_j$.  In view of Observation 2.3, minimality of the pair $(\udim T_j, \, \udim \soc M_j)$ in $\topsocbd$ guarantees moreover that $\soc M_j$ is the generic socle of  $\C_j$.   In particular, $\C_i = \C_j$ implies $i = j$.   
\qed \enddemo

By Corollary 3.7, finding the irreducible components of $\modlad$ amounts to finding the  minimal pairs in $\topsocbd$.

\definition{Example 3.5 revisited} Let $\la$ be as in 3.5. Inspection of the set $\topsocbd$ yields the following distinct irreducible components of $\modlad$ for $\bd = d = 13$:  Namely the closures $\overline{\toptd}$, where $T= S^t$ with $4 \le t \le 9$.  Theorem 3.12 below will place this example into a general context.
\enddefinition

\proclaim{Corollary 3.8}  Suppose that $T \in \lamod$ is semisimple.  Whenever there exists a module $M$ with dimension vector $\bd$ and top $T$ which has no simple direct summand, the closure of $\toptd$ in $\modlad$ is an irreducible component of $\modlad$. \endproclaim

\demo{Proof}  Indeed, if there exists a module $M$ as specified, $\udim \top M + \udim \soc M = \bd$, which guarantees that the pair $(\udim \top M, \, \udim \soc M)$ is minimal in $\topsocbd$. \qed \enddemo 

While Corollary 3.7 provides us with an algorithm to compute the irreducible components of the varieties $\modlad$,  the road by way of a full calculation of the set $\topsocbd$ is unnecessarily labor-intensive.  In a nutshell, the upcoming, far more convenient, approach to the components may be paraphrased as follows:  A semisimple module $T$ is a generic top of $\modlad$ precisely when no simple summand $S$ of $T$ may be shifted into the radical of a module with top $T/S$.  

For a given dimension vector $\bd$, generic modules with, resp\. without, simple direct summands will usually coexist in $\modlad$.  Here is an illustrative example of low dimension:  Let 
$$\la = KQ/\langle \text{paths of length\ } 2\rangle,$$ 
where $Q$ is the quiver 
$$\xymatrixrowsep{2pc}\xymatrixcolsep{2pc}
\xymatrix{
1 \ar[r] \ar@/_1pc/[rr] &2 \ar[r] &3
}$$
\noindent and let $\bd = (1, 1, 1)$.  Then $\modlad$ has two irreducible components, namely the closures of the orbits of the modules $M_1 = \la e_1$ and $M_2 = \la e_2 \oplus S_1$.  The corresponding minimal pairs in $\topsocbd$ are $\bigl( (1, 0, 0), (0, 1, 1) \bigr)$ and $\bigl( (0, 1, 1), (1, 0, 1) \bigr)$.  If, on the other hand, we introduce an additional arrow from $2$ to $1$ into $Q$, the two generic modules of $\modlad$ are the indecomposables $\la e_1$ and $\la e_2$.

\proclaim{Proposition 3.9.  Finding the generic tops of $\modlad$}  We continue to assume that $J^2 = 0$.  For a dimension vector $\bd$ and a semisimple module $T$, the following statements are equivalent:
\smallskip

{\rm \bf{(a)}} $T$ is a generic top of $\,\modlad$, i.e., $\overline{\toptd}$ is an irreducible component of $\modlad$.   
\smallskip

{\rm \bf{(b)}}  There exists a $\la$-module $M$ with top $T$ and dimension vector $\bd$ having the following property:  Whenever $S \in \lamod$ is a simple direct summand of $M$, say $M = M' \oplus S$, and $P' = P(M')$ is a projective cover of $M'$, the semisimple module $JM \oplus S$ fails to embed into $JP'$.
\endproclaim

\demo{Proof}  To prove ``(a)$\implies$(b)", suppose that $T$ is a generic top of $\modlad$.  Corollary 3.7 then yields a module $M$ with $\top M = T$ such that the pair $(\udim \top M,\, \udim \soc M)$ is minimal in $\topsocbd$.  To show that $M$ has the property specified under (b), let $S$ be a simple direct summand of $M$, say $M = M' \oplus S$, and $P'$ a projective cover of $M'$.  Clearly, $JM = JM'$.  We assume that $JM' \oplus S \subseteq J P'$ and let $C'$ be any $\la$-direct complement of $JM' \oplus S$ in the semisimple module $JP'$.   Then the quotient $N = P'/C'$ in turn has dimension vector $\bd$.  However, $\top N$ is properly contained in $\top M$, while $\soc N = \soc M' \oplus S = \soc M$.  Therefore $(\udim \top N,\, \udim \soc N) < (\udim \top M,\, \udim \soc M)$, which contradicts our choice of $M$.

Now suppose that $M$ is a module satisfying condition (b).  To deduce (a), we will show that $\top M = \top N$ for some module $N$ with the property that $(\udim \top N, \, \udim \soc N)$ is minimal in $\topsocbd$, and then apply Corollary 3.7.  Indeed, choose a minimal element $(\udim \top N,\, \udim \soc N)$ in $\topsocbd$, subject to the inequality 
$$(\udim \top N,\, \udim \soc N) \le  (\udim \top M,\, \udim \soc M),$$ 
and let $U$ be semisimple such that $\top N \oplus U = \top M$.  In light of  $\udim M = \udim N$, we deduce $J N = J M \oplus U \subseteq \soc N$ $\subseteq$ $ \soc M$.  Consequently, $U$ is (isomorphic to) a direct summand of $M$, say $M \cong \Mhat \oplus U$.  Let $\Phat$ be a projective cover of $\Mhat$, which also makes $\Phat$ a projective cover of $N$.  In light of our hypothesis on $M$, the inclusion $JM \oplus U = JN \subseteq J \Phat$ then forces $U$ to be zero.  In other words, we conclude $\top M = \top N$ as desired.
\qed   
\enddemo 

Note that Corollary 3.8 is in turn subsumed in Proposition 3.9 as an obvious special case. The easy proofs of the following remarks  --  helpful in applying Proposition 3.9  --   are left to the reader.

\definition{Supplement 3.10}  Suppose $T= \bigoplus_{1 \le i \le n} S_i^{t_i}$ is  realizable with respect to $\bd$.  Let $U$ be the generic socle of the modules in $\toptd$, and $E$ the injective envelope of $\bigoplus_{1 \le i \le n} S_i^{d_i - t_i}$.  

\noindent {\bf{(a)}}  Generically, the modules in $\toptd$ have a direct summand isomorphic to $S_k$ if and only if $S_k^{t_k}$ fails to be contained in the factor module $E / \soc E$.  In particular:  The modules in $\toptd$ are generically without simple direct summands if and only if  $T \subseteq E/\soc E$ (cf\. Observation 3.4).  
\smallskip

\noindent {\bf{(b)}}  Let $D$ be the standard duality $\hom_K( - , K): \lamod \rightarrow \text{mod-}\la$, and let $\bd$ also stand for the dimension vector of the right $\la$-module $D \bigl(\bigoplus_{1 \le i \le n}  S_i^{d_i} \bigr)$. 

Then $T$ is a generic top of $\modlad = \rep_{\bd}(_{\la}\la)$ if and only if $D(U)$ is a generic top of $\rep_{\bd}(\la_{\la})$.  Moreover, $D(T)$ is the generic socle of the right $\la$-modules with dimension vector $\bd$ and top $D(U)$ (cf\. Corollary 3.7). \qed
\enddefinition

We follow with an example to illustrate the extra computational edge we gain from Proposition 3.9.          

\definition{Example 3.11}  Let $\la = KQ/\langle \text{paths of length\ } 2\rangle$, where $Q$ is the quiver
$$\xymatrixrowsep{4pc}\xymatrixcolsep{4pc}
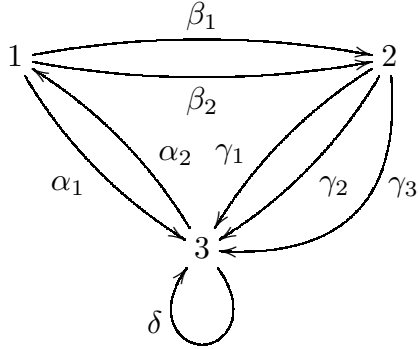
\xymatrix{
1 \ar@/_/[dr]_{\alpha_1} \ar@/^/[rr]^{\beta_1} \ar@/_/[rr]_{\beta_2} &&2 \ar@/_/[dl]_(0.65){\gamma_1} \ar@/^/[dl]^{\gamma_2} \ar@/^2.2pc/[dl]^(0.35){\gamma_3}  \\
 &3 \ar@'{@+{[0,0]+(6,-6)} @+{[0,0]+(0,-15)} @+{[0,0]+(-6,-6)}}^(0.75){\delta} \ar@/_/[ul]_(0.35){\alpha_2}
}$$
Using Proposition 3.9, one finds that, for $\bd = (3,3,3)$, the variety $\modlad$ has precisely $16$ irreducible components with generic tops as follows:  $S_2^3 \oplus S_3^3$, $S_1 \oplus S_2 \oplus S_3^2$, $S_1 \oplus S_2^2 \oplus S_3^2$, $S_1 \oplus S_2^3 \oplus S_3^2$, $S_1^2 \oplus S_2 \oplus S_3$, $S_1^2 \oplus S_2 \oplus S_3^2$, $S_1^2 \oplus S_2^2 \oplus S_3$, $S_1^2 \oplus S_2^2 \oplus S_3^2$, $S_1^2 \oplus S_2^3 \oplus S_3$, $S_1^3$, $S_1^3 \oplus S_3$, $S_1^3 \oplus S_2$, $S_1^3 \oplus S_2 \oplus S_3$, $S_1^3 \oplus S_2^2$, $S_1^3 \oplus S_2^2 \oplus S_3$, $S_1^3 \oplus S_2^3$.

First, we indicate how to prune the set of semisimple modules $S_1^i \oplus S_2^j \oplus S_3^k$ with $0 \le i,j,k \le 3$ by discarding those that fail to be realizable with respect to $\bd$.  Then we give a few sample arguments for the decision process, whether a given realizable semisimple module is a generic top of $\modlad$.  One readily extrapolates to obtain an algorithm.

For $i = 0$, i.e., $T = S_2^j \oplus S_3^k$, the projective cover $P$ of $T$ contains the simple $S_1$ with multiplicity $k$ and the simple $S_2$ with multiplicity $j$.  Thus $T = S_2^3 \oplus S_3^3$ is the only realizable choice for $i = 0$.  Proposition 3.9 shows $T$ to be, in fact, a generic top of $\modlad$:  Indeed, $M = S_2^3 \oplus (\la e_3/ \la \delta)^3$ has top $T$ and dimension vector $\bd$; whenever $S \subseteq M$ is one of the direct summands isomorphic to $S_2$, we denote by $P'$ the projective cover of $M/S$ and find $S \not\subseteq JP'$; a fortiori $JM \oplus S \not\subseteq JP'$.  

For $i = 1$, realizabilty of $S_1 \oplus S_2^j \oplus S_3^k$ entails $j \ge 1$ and $k \ge 2$.  Let us consider the smallest choice, $T =  S_1 \oplus S_2\oplus S_3^2$.  Since the injective envelope $E$ of $S_1^2 \oplus S_2^2 \oplus S_3$ has the property that $T \subseteq  E / \soc E$, the modules in $\toptd$ are generically devoid of simple direct summands (see Supplement 3.10), and again Proposition 3.9 guarantees that $T$ is a generic top of $\modlad$.  On the other hand, increasing the top to $T = S_1 \oplus S_2 \oplus S_3^3$, we obtain a realizable semisimple, which fails to be a generic top of $\modlad$:  Indeed, any module $M$ with top $T$ and dimension vector $\bd$ has radical $JM = S_1^2 \oplus S_2^2$, and the injective envelope $E$ of $JM$ contains $S_3$ only with multiplicity $2$ in its top.  Thus, $M$ has a direct summand $S_3$, say $M \cong M' \oplus S_3$; however, the radical of the projective cover of $M'$ {\it does\/} contain $JM \oplus S_3$.  

Yet, there are generic tops of $\modlad$ which have dimension larger than $5$, e.g., $T = S_1^2 \oplus S_2^3 \oplus S_3$. For verification, apply Corollary 3.8 to the module with graph
$$\xymatrixrowsep{1pc}\xymatrixcolsep{0.75pc}
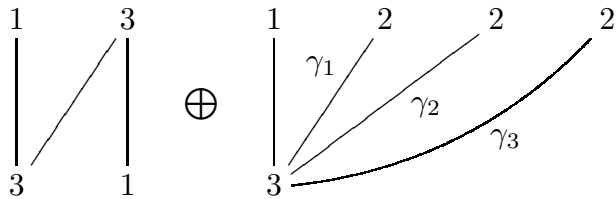
\xymatrix{
1 \edge[dd] &&3 \edge[ddll] \edge[dd] &&1 \edge[dd] &&2 \edge[ddll]_(0.35){\gamma_1} &&2 \edge[ddllll]^(0.4){\gamma_2} &&2 \edge@/^1pc/[ddllllll]^(0.4){\gamma_3}  \\
 &&&\bigoplus  \\
3 &&1 &&3
 }$$
\enddefinition

\noindent {\bf{(C)  Local algebras with $J^2 = 0$}} 
\smallskip

Now suppose that $\la$ is local, meaning that the quiver $Q$ consists of a single vertex and $r$ loops:
$$\xymatrixrowsep{0.75pc}\xymatrixcolsep{0.3pc}
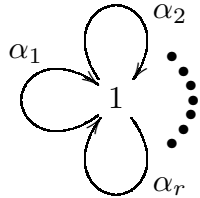
\xymatrix{
 & \ar@{{}{*}}@/^1pc/[dd] \\
1 \ar@'{@+{[0,0]+(-6,-6)} @+{[0,0]+(-15,0)} @+{[0,0]+(-6,6)}}^(0.7){\alpha_1}
\ar@'{@+{[0,0]+(-6,6)} @+{[0,0]+(0,15)} @+{[0,0]+(6,6)}}^(0.7){\alpha_2}
\ar@'{@+{[0,0]+(6,-6)} @+{[0,0]+(0,-15)} @+{[0,0]+(-6,-6)}}^(0.3){\alpha_r}  \\
 &
}$$
For $J^2 = 0$ and any positive integer $d$, the variety $\Mod_d(\la)$ then consists of the sequences $(A_1, \dots, A_r)$ of linear operators $A_i \in \End_K(K^d)$ with the property that $A_i A_j = 0$ for all $i,j$.  In this situation, the irreducible components of $\Mod_d(\la)$ can be described in a unified format, in terms of $r$ and $d$, which yields a closed formula for the number of components.  We preempt the general discussion with that of the trivial case, $r =1$, which plays an outsider role.  In this case, i.e., when $\la \cong K[X]/ (X^2)$, each of the $\Mod_d(\la)$ is irreducible.  Indeed, if $d = 2k$ is even, then $\Mod_d(\la)$ is the orbit closure of the free left $\la$-module $\la^k$; for $d = 2k+1$, $\Mod_d(\la)$ is the orbit closure of $\la^k \oplus \bigl( \la/(X) \bigr)$.  
  
\proclaim{Theorem 3.12}  Suppose $\la$ is a basic local finite dimensional $K$-algebra with $J^2 = 0$, say $\la = KQ/\langle \text{all paths of length}\ 2 \rangle$, where $Q$ is the quiver with a single vertex and $r \ge 2$ distinct loops.  Let $S$ be the unique simple in $\lamod$, and, for $d \ge 2$, let $u \in \NN$ be minimal with respect to $r \cdot u \ge d-u$, i.e., $u$ is the integer ceiling of $d / (r+1)$.  

Then the irreducible components of $\Mod_d(\la)$ are precisely the closures $\overline{\rep^T_d}$, where $T= S^t$ with $u \le t \le d - u$.  Moreover, the modules in the irreducible components of $\Mod_d(\la)$ are generically without simple direct summands. 
\endproclaim 

\demo{Proof}  We start by observing that the pairs $(t, d- t)$, where $t \in \NN \cap [u, d-u]$ are precisely those pairs $(a, b) \in \NN \times \NN$ with $a+
b = d$ such that $a \le r \cdot b$ and $b \le r \cdot a$.   In the following, let $P$ be a projective cover of $S^t$ and $E$ an injective envelope of $S^{d-t}$.

First suppose that $T = S^t$ with $t \in [u, d - u]$.  Since $\dim JP = r \cdot t \ge d-t$, we find that $T$ is realizable with respect to $d$.  Moreover, from $\dim (E/ \soc E) = r \cdot (d- t)$, we infer that $E$ contains a submodule $M$ with $\SS(M) = (S^t, S^{d - t})$.  This shows that, generically, the modules with top $T$ embed into $E$, which in turn implies that, generically, they are free of simple direct summands.  Therefore, $T$ is a generic top of $\Mod_d(\la)$ by Corollary 3.8.

Now suppose that $T = S^t$ with $t \notin [u,\, d-u]$.  If $t < u$, then $T$ fails to be realizable with respect to $d$.  If, on the other hand, $t > d-u$, then any $d$-dimensional $\la$-module $M$ with top $T$ satisfies $\dim \top M > r \cdot \dim JM$, whence $M$ does not embed into $E$.  Consequently $M$ has a simple direct summand.  On the other hand, clearly $\dim JP \ge \dim JM + 1$, whence $JM \oplus S$ embeds into $JP$.  Applying Proposition 3.9, we thus conclude that $T$ fails to be a generic top of $\Mod_d(\la)$.
\qed \enddemo

Letting $\la_r$ be the basic local $K$-algebra with $\dim J = r$, we deduce that, for any dimension $d \ge 2$, the number of irreducible components of $\Mod_d(\la_r)$ asymptotically grows like $d - 1$ for $r \rightarrow \infty$.  More precisely, we obtain: 

\proclaim{Corollary 3.13}  {\rm{(See \cite{\DoFl} and \cite{\Mor, Section 5} for the case $r=2$.)}}  Again suppose that $\la = KQ/\langle \text{all paths of length}\ 2 \rangle$, where $Q$ has a single vertex and $r \ge 2$ distinct loops.  For any integer $d \ge 2$, the number of irreducible components of $\Mod_d(\la)$ is
$$d \ - \ 2 \bigg \lceil \frac{d}{r+1} \bigg \rceil \ +\  1.$$
\endproclaim

\demo{Proof}  This is an immediate consequence of Theorem 3.12, since the displayed value just counts the natural numbers in the interval $[u, \, d-u]$. \qed \enddemo 

\definition{Example 3.14}  For $r = 3$ and $d = 8$, the generic radical layerings which bijectively tag the five irreducible components of $\Mod_d(\la)$ may be visualized as follows.  Here each bullet stands for a composition factor $S$.
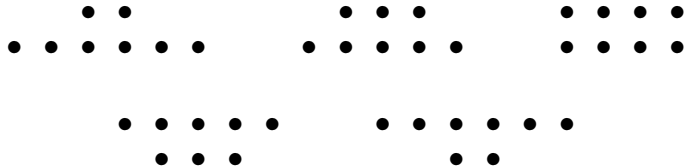
$$\xymatrixrowsep{0.05pc}\xymatrixcolsep{0.05pc}
\xymatrix{
 &&\bullet &\bullet && &&& &\bullet &\bullet &\bullet & &&&\bullet &\bullet &\bullet &\bullet  \\
\bullet &\bullet &\bullet &\bullet &\bullet &\bullet &&&\bullet &\bullet &\bullet &\bullet &\bullet &&&\bullet &\bullet &\bullet &\bullet  \\  \\  \\
 &&&\bullet &\bullet &\bullet &\bullet &\bullet &&&\bullet &\bullet &\bullet &\bullet &\bullet &\bullet  \\
 &&& &\bullet &\bullet &\bullet & &&& &&\bullet &\bullet
 }$$
\enddefinition

\head 4.  The modules parametrized by the irreducible components of $\modlad$.  Geometry of the components \endhead

We continue to restrict to basic finite dimensional algebras with $J^2 = 0$.  In this brief section, we explicitly describe the modules parametrized by the closure of $\toptd$ in $\modlad$ for any semisimple $T$.  We then use this information to describe downsized and simplified affine and projective parametrizing varieties for these closures, apt to facilitate further investigation into the representations parametrized by the irreducible components (see, for instance, Section 6 below).  
\bigskip

\noindent {\bf (A)  Recognizing the modules in $\overline{\toptd}$} 
\smallskip

Since all irreducible components of the varieties $\modlad$ are of the form $\overline{\toptd}$ for suitable choices of $T$, Corollary 4.2 below yields, in particular, a description of the modules in the irreducible components of $\modlad$, once the generic tops have been determined.     

\proclaim{Proposition 4.1}  We retain the hypothesis $J^2 = 0$.  Moreover, we let $T \in \lamod$ be a semisimple module which is realizable with respect to the dimension vector $\bd$.  Then:
\smallskip

{\rm \bf{(a)}}  Every module $X$ in $\overline{\toptd}$ is a degeneration of some module  in $\toptd$.
\smallskip

{\rm \bf{(b)}}  Let $X \in \lamod$ with $T \subsetneqq \top X$.  Then $X$ is a degeneration of some module in $\toptd$ if and only if $X$ decomposes in the form
$$X \cong X' \oplus U, \tag\ddagger$$
where $X'$ has top $T$ and $U$ is a nonzero semisimple module.
\endproclaim

\demo{Proof}  We will first show that every module $X$ which belongs to $\overline{\toptd} \setminus \toptd$ is of the form ($\ddagger$) displayed in part (b). 

So let $X$ be in the closure of $\toptd$, with $T \subsetneqq \top X$.  Say $\top X = T \oplus U$.  Moreover, let $N$ be any module in $\toptd$ whose socle is the generic one in $\toptd$ (and thus also in $\overline{\toptd}$).  Write $N = N_0 \oplus N_1$, where $N_1$ is semisimple and $N_0$ is without simple direct summands.    From the equality $\udim X = \udim N$ we deduce $JN = JX \oplus U$, which entails 
$$JX \oplus U \oplus N_1 =  JN \oplus N_1 = JN_0 \oplus N_1 = \soc N \subseteq \soc X;$$
the equalities follow from Observation 3.4, and Observation 2.3 yields the final inclusion.  We infer that $X$ contains a copy $V$ of the semisimple module $U \oplus N_1$ which meets $JX$ trivially and conclude that $V$ is a direct summand of $X$, say $X \cong X'' \oplus V$.   Consequently, $X \cong X' \oplus U$, where $X' = X'' \oplus N_1$ has top $T$ by construction.  In particular, this shows that every non-top-stable degeneration $X$ of some module in $\toptd$ has the form postulated in (b).    

To complete the proofs of (a) and (b), it now suffices to check that every module $X = X' \oplus U$ as specified in $(\ddagger)$ is a degeneration of some module $M$ in $\toptd$.  Let $P$ be a projective cover of $T$ and hence also of $X'$; say $X' = P/D$ for some submodule $D \subseteq JP$.  Realizability of $T$ with respect to $\bd$ implies that $D \oplus U$ embeds into $JP$.  Choose a $\la$-direct complement $C$ of $D \oplus U$ in $JP$, and note that $X$ is a degeneration of the module $M = P/C$, the latter being an object of $\toptd$.   
\qed \enddemo

We derive the following simplistic description of the modules in the closure of $\toptd$.

\proclaim{Corollary 4.2}  Suppose $J^2 = 0$.  Let $T$ be a semisimple $\la$-module which is realizable with respect to $\bd$.  Then the modules in the closed subvariety $\overline{\toptd}$ of $\modlad$ are precisely those modules with dimension vector $\bd$ which have the form $M \oplus U$, where $M$ has top $T$ and $U$ is semisimple. \qed
\endproclaim

\bigskip

\noindent {\bf (B)  Alternative affine parametrization of the modules in $\overline{\toptd}$} 
\smallskip

Once again, the upcoming results apply, in particular, to the irreducible components of the varieties $\modlad$.  

We start by decomposing $T$ into its homogeneous components: $T = \bigoplus_{1 \le i \le n} T_i$, where $T_i = S_i^{t_i}$, and let $R = \bigoplus_{1 \le i \le n} R_i$, where $R_i = S_i^{d_i - t_i}$.  Then the closure of $\toptd$ in $\modlad$ may be parametrized by the affine space
$$\C(T, \bd) = \{(x_\alpha)_{\alpha \in Q_1} \mid x_\alpha \in \hom_K (T_{\start(\alpha)}, R_{\term(\alpha)})\},$$
via the assignment  
$$\rho: x\   \mapsto \   {}_\la(T \oplus R), \   \text{where}\ \, \alpha(v + w) = x_\alpha (v_{\start(\alpha)}) \ \ \text{for}\  \ v = v_1 + \cdots + v_n \in T,\,  w \in R.$$ 
Next we focus on the map, induced by this parametrization, from $\C(T,\bd)$ to the set of isomorphism classes of $\la$-modules with dimension vector $\bd$.  As is evidenced by Corollary 4.2, the image of this map consists precisely of the isomorphism classes of modules in $\overline{\toptd}$.  

In addition, we consider the reductive group 
$$G(T) \ = \ \biggl(\prod_{1 \le i \le n} \GL(T_i)\biggr) \, \times \biggl( \prod_{1 \le i \le n} \GL(R_i)\biggr),$$ 
and let it act on $\C(T,\bd)$ via 
$$\bigl((g_i), (h_i) \bigr) . x \ = \bigl(h_{\term(\alpha)} \, x_\alpha \,  g^{-1}_{\start(\alpha)}\bigr).$$ 
The original affine parametrization of $\overline{\toptd}$ can thus be whittled down to a parametrization by a comparatively small affine space.  We record this fact as

\proclaim{Observation 4.3}  The image of the map 
$$\rho: \C(T,\bd) \rightarrow \{\text{isomorphism classes of modules with dimension vector}\ \bd\}$$
equals the set of isomorphism classes of modules in $\overline{\toptd}$, and the fibers of $\rho$ coincide with the orbits of the action of $G(T)$ on $\C(T, \bd)$. \qed
\endproclaim

It is well-known that the algebras with vanishing radical square are stably equivalent to hereditary algebras.  However, while for any hereditary algebra $\lahat$ the module varieties $\Mod_{\bd}(\lahat)$ are affine spaces, this is far from being true for algebras $\la$ with $J^2 = 0$:  Not only do we know the $\modlad$ to admit arbitrarily high numbers of irreducible components, the individual components may have singularities (see, e.g., \cite{\GoHu, last page}).  On the other hand, the strong homological ties between the two classes of algebras {\it are\/}, in fact, paralleled by geometric ones:  The first piece of evidence, Observation 4.3, is only the tip of the iceberg.  The geometric connection will become more salient in Section 5, where we focus entirely on the projective varieties parametrizing modules with fixed top.  (In Subsection 4.C, the links to the hereditary case will not surface yet.) 
\bigskip

\noindent {\bf (C)  Alternative projective parametrization.  Classification}} 
\smallskip

Once again, we point out that the following results target, in particular, the irreducible components $\overline{\biggrasstd}$ of $\biggrasslad$, where $T$ traces the generic tops of the modules with dimension vector $\bd$.  

We still let $T$ be any semisimple module which is realizable with respect to $\bd$ and, as usual, we identify isomorphic semisimple modules.  To take full advantage of the increased transparency of the ``small" Grassmannians $\grasstbd$ parametrizing the modules with top $T$, we subdivide the closure of $\biggrasstd$ in $\biggrasslad$ into portions, each of which encodes the objects $M \oplus U$ for a fixed semisimple summand $U$ of $\bigoplus_{1 \le i \le n} S_i^{d_i - t_i}$ (in the notation of Corollary 4.2, i.e.,  $\top M = T$).  This stratification of the closure also advertises itself from the viewpoint of a potential classification of the modules in $\overline{\biggrasstd}$, in a sense made precise below.  In fact, Observation 4.6 shows classifiability of the isomorphism classes of modules in the individual segments to be ``the" best one can hope for by way of a geometric approach.  

For every submodule $U \subseteq \bigoplus_{1 \le i \le n} S_i^{d_i - t_i}$, we consider the set of isomorphism classes
$$\M(\bd, T, U) = \{ [M \oplus U] \mid  \udim M = \bd - \udim U,\, \top M = T\};$$ 
here $[X]$ stands for the isomorphism class of a module $X$.  Note that the special case $U= 0$ takes us back to the modules in $\biggrasstd$ or, equivalently, in $\grasstd$.  By Corollary 4.2, the disjoint union of the $\M(\bd, T, U)$ equals the full set of isomorphism classes of modules represented by $\overline{\biggrasstd}$.  Clearly, each $\M(\bd, T, U)$ is parametrized by the projective variety $\grass^T_{(\bd - \udim U)}$, via 
$$C \mapsto \text{the class of}\  (P/C) \oplus U, \ \ \ \text{for}\ \  C \in \grass^T_{(\bd - \udim U)},$$
where $P$ is the distinguished projective cover of $T$ (cf\. Section 2)

The parametrizing varieties $\grass^T_{(\bd - \udim U)}$ of the individual portions $\M(\bd, T, U)$ are comparatively simplistic from a geometric perspective.  Namely, they are all direct products of classical Grassmann varieties:

\proclaim{Observation 4.4.  Structure of $\grasstbd$}  Again assume $J^2 = 0$.  Suppose $T = \bigoplus_{1 \le i \le n} S_i^{t_i}$ is realizable with respect to $\bd$, and let $P$ be the projective cover of $T$.  If $JP = \bigoplus_{1 \le i \le n} S_i^{p_i}$, then $\grasstbd$ is the following direct product of classical vector space Grassmannians:
$$\grasstbd\ \ \   \cong \ \ \ \prod_{1 \le i \le n} \Gr (d_i -  t_i, K^{p_i}). \ \ \ \qed$$
\endproclaim

The following is an immediate consequence of \cite{\classifying, Corollary 4.5}.  For background on fine/coarse moduli spaces, we refer to \cite{\moduli}.

\proclaim{Observation 4.5. Classification of the modules in $\overline{\biggrasstd}$ when $T$ is simple}  Suppose $J^2 = 0$.  If $T$ is a simple $\la$-module which is realizable with respect to $\bd$, then each of the sets $\M(\bd, T, U)$ has a fine moduli space, namely the variety $\grass^T_{(\bd - \udim U)}$.  The corresponding universal family is {\rm{(}}informally presented{\rm{)}}
$$ \bigl( P/C \oplus U \bigr), \ \ \text{where \ } C \text{\ traces\ } \grass^T_{(\bd - \udim U)}.$$

In particular, the modules in an irreducible component $\C$ of $\biggrasslad$ are classifiable by a finite number of projective moduli spaces whenever $\C$  has simple generic top.
 \qed \endproclaim 
 
 \proclaim{Observation 4.6.  Limitations to classifiability by a single moduli space}  
Still $J^2 =  0$.  Let  $T = \bigoplus_{1 \le i \le n} S_i^{t_i}$ be realizable with respect to $\bd$, take $U \subseteq \bigoplus_{1 \le i \le n} S_i^{d_i - t_i}$, and suppose $\M$ is a set of isomorphism classes represented by a subvariety of $\biggrasslad$, with the property that  $\M(\bd, T, U) \subseteq \M$.  

If the objects in $\M$ possess a coarse moduli space classifying them up to isomorphism, then $\M$ does not intersect any $\M(\bd, T, U')$ with $U \subsetneqq U'$ or $U' \subsetneqq U$. \endproclaim 

\demo{Proof}  Let $\X$ be the union of those $\aut_\la(\bP)$-orbits in $\biggrasslad$ which correspond to the isomorphism classes in $\M$, and $\X(\bd, T, U)$ the union of the orbits corresponding to the classes in $\M(\bd, T, U)$.  If $\M$ is classifiable by a coarse moduli space, then all orbits contained in $\X$ are relatively closed in $\X$.  First suppose that $U \subsetneqq U' \subseteq \bigoplus_{1 \le i \le n} S_i^{d_i - t_i}$.  Since every module in $\M(\bd, T, U')$ is a degeneration of some module in $\M(\bd, T, U)$  --  this follows from Proposition 4.1, with $T$ is replaced by $T \oplus U$  --   we deduce from the inclusion $\X(\bd, T, U) \subseteq \X$ that $\X \cap \X(\bd, T, U')  = \varnothing$.  Next suppose that $U' \subsetneqq U \subseteq \bigoplus_{1 \le i \le n} S_i^{d_i - t_i}$.  Since every module in $\M(\bd, T, U')$ degenerates to one in $\M(\bd, T, U)$,  we arrive at the same conclusion.  \qed \enddemo

\head 5.  Relating the module varieties of $\la$ to those of a stably equivalent hereditary algebra \endhead

\noindent 
Throughout this section we assume that $\Lambda$ is an algebra with vanishing radical
square.  In other words,
$$\la = KQ/ I \ \ \text{where} \ \ I=\langle\text{all paths of length\ }2\rangle \ \ \text{for some quiver \ } Q.$$
In Section 3, we saw how to determine the tops $T$ which are generic with respect to any given dimension vector $\bd$, which allows us to identify the irreducible components $\overline{\biggrasstd}$ of the variety $\biggrasslad$. (Recall that this is equivalent to identifying the components of $\modlad$.)  The next step in our program is to explore generic data regarding the $\la$-modules encoded by the components. Clearly, this amounts to assembling generic information on the modules in $\biggrasstd$ or, in other words, assembling generic information on the modules in $\grasstbd$.  It will turn out that this second task may be tackled by tapping into the solidly developed theory of hereditary algebras.   

In fact, we will establish a geometric counterpart to the well-known fact that any algebra with vanishing radical square is stably equivalent to a hereditary algebra $\lahat$.  There is a network of bridges, each connecting varieties which parametrize modules with fixed dimension vector and fixed top over the algebras we are comparing.  These bridges thus provide exactly the logistics required for our purpose.  

Following standard practice, we choose $\lahat$ to be $K\Qhat$, where $\Qhat$ is the {\it separated quiver of $Q = (Q_0, Q_1)$\/}; see [\AuReSm, p\. 350].  Recall the definition: $\Qhat = (\Qhat_0, \Qhat_1)$, where $\Qhat_0$ is the disjoint union of $Q_0 = \{e_1, \dots, e_n\}$ and a duplicate of $Q_0$, written as $\{\ehat_1, \dots, \ehat_n\}$, so that the vertex set $\Qhat_0$ has cardinality $2n$. The set $\Qhat_1$ of arrows duplicates $Q_1$; it is written in the form $\{\alphahat \mid \alpha \in Q_1\}$, where $\alphahat$ is an arrow from the vertex $e_i$ to the vertex $\ehat_j$ of $\Qhat$, provided that $\alpha \in Q_1$ is an arrow from $e_i$ to $e_j$.   

Clearly, the vertices $\ehat_j$ are all sinks of $\Qhat$, whence $\Qhat$ has no paths of length larger than $1$.  In particular, $\lahat$ is a finite dimensional hereditary algebra whose radical $\Jhat$ in turn has vanishing square.  Note that the indecomposable projective modules $\lahat \ehat_i$, for $1 \le i \le n$, are simple. 
We will continue to write the simples in $\lamod$ as $S_i$; the simples in $\lahatmod$ which correspond to the ``old" vertices $e_1, \dots, e_n$ are denoted by $\Shat(e_i) = \lahat e_i/ \Jhat e_i$, those corresponding to the ``new" vertices $\ehat_1, \dots, \ehat_n$, by $\Shat(\ehat_j) = \lahat \ehat_j$.  In order to obtain dimension vectors of $\lahat$ which transparently relate to those of $\la$, we order the vertices in $\Qhat$ (and accordingly the simple $\lahat$-modules) as follows: $e_1, \dots, e_n, \ehat_1, \dots, \ehat_n$.

Instead of using the standard triangular matrix functor, employed in \cite{\AuReSm} to show that the categories $\lamod$ and $\lahatmod$ are indeed stably equivalent, it will for our purposes be preferable to specify correspondences between sets of isomorphism classes of $\la$- resp\. $\lahat$-modules set apart by their projective presentations.  These correspondences will parallel our description of isomorphisms linking suitable pairs of parameter varieties for $\la$- and $\lahat$-modules.

The $\lahat$-modules of interest are those induced from $\la$ in the following sense.

\definition{Definition 5.1 of induced modules}    A $\lahat$-module $N$ will be called {\it induced from $\la$\/} (or simply {\it induced\/}) in case the top $N/ \Jhat N$ is a direct sum of copies of the simples $\Shat(e_i)$, $1 \le i \le n$.  
\enddefinition 

The following remarks are obvious:

\proclaim{Lemma 5.2}  $\bullet$ For every $\lahat$-module $N$, the radical $\Jhat N$ is a direct sum of copies of the $\Shat(\ehat_j)$.

$\bullet$ Every $\lahat$-module $N$ is {\rm{(}}uniquely{\rm{)}} a direct sum of a module induced from $\la$ and a direct sum of copies of projective simples $\Shat(\ehat_j)$. In particular, every non-simple indecomposable $\lahat$-module is induced from $\la$.  \qed
\endproclaim

\medskip

\noindent {\bf{(A)  A two-way shift of geometric information $\lamod$ $\leftrightarrow$ $\lahatmod$. }} 
\smallskip

We begin with a precise description of the matchup between (isomorphism classes of) $\la$-modules and (isomorphism classes of) induced $\lahat$-modules.  To that end, we introduce some further notation, guided by the intuitive picture.  We start with a semisimple $\la$-module $T$ with dimension vector $(t_1, \dots, t_n)$ of total dimension $t$.  As in the definition of the projective variety $\grasstd$ in Section 2, we fix a projective cover $P = \bigoplus_{1 \le r \le t} \la z_r$ of $T$, where $z_r = e(r) z_r$ are top elements of $P$.  We will write $z_r = \e(r)$, so as to emphasize the norming idempotent $e(r) \in \{e_1, \dots, e_n\}$.  Moreover, we consider the twin projective $\lahat$-module $\Phat = \bigoplus_{1 \le r \le t} \lahat \e(r)$ with top elements $\e(r)$ normed by the analogous idempotents $e(r)$, now viewed as elements of $\lahat$.   It is obvious that the {\it nonzero\/} elements $\alpha \e(r) \in P$, where $\alpha$ traces the ($I$-residue classes of) arrows in $Q_1$, form a $K$-basis for $JP$, the analogous statement being true for the nonzero $\alphahat \e(r) \in \Jhat \Phat$.  Note that 
$\alpha \e(r) \ne 0$ precisely when $\alphahat \e(r) \ne 0$. 

For any element $c \in JP$, that is, for any $K$-linear combination $c$ of such basis elements $\alpha \e(r)$ of $JP$, we let $\chat$ be the corresponding $K$-linear combination of the elements $\alphahat \e(r)$ in $\Jhat \Phat$.  In case $M = P/C$, where $C$ is a submodule of $JP$, we define $\Chat$ to be the $\lahat$-submodule  of $\Jhat \Phat$ consisting of the elements $\chat$ for $c\in C$, and set $\Mhat = \Phat/\Chat$.  Evidently, the $\lahat$-module $\Mhat$ is then induced from $\la$. For the semisimple module $T = \bigoplus_{1 \le i \le n} S_i^{t_i}$, the corresponding induced module $\That =  \bigoplus_{1 \le i \le n} \Shat(e_i)^{t_i}$ is in turn semisimple.

We thus obtain a well-defined assignment of isomorphism classes 
$$M = P/C \mapsto \Mhat = \Phat/\Chat.$$ 
Our setup entails that any $\lahat$-module $\Mhat$ which is induced from $\la$ is isomorphic to $\Phat/ \widehat{C}$ for some projective $\la$-module $P$ and submodule $C \subseteq JP$.  As we will ascertain, the above assignment gives rise to a well-defined bijection between the set of isomorphism classes of $\la$-modules with top $T = P/JP$ on one hand and the set of isomorphism classes of $\lahat$-modules with top $\That = \Phat/ \Jhat \Phat$ on the other.  This pairing of isomorphism types of modules is paralleled by a family of isomorphisms connecting the pertinent parametrizing varieties, the latter isomorphisms well behaved relative to the acting groups, $\autlap$ and $\autlaphat$.  Clearly, $\autlap$ will have higher dimension than $\autlaphat$ in general, but the difference is erased as far as the effect of the action is concerned.  For more precision, see Proposition 5.3 below.  

It will be convenient to identify the factor group $\autlap/(\autlap)_u$ of $\autlap$ modulo its unipotent radical with $\aut_\la(P/JP) = \autlat$; this is harmless in light of the fact that  $\autlap \cong \autlat \ltimes (\autlap)_u$.  The  automorphism group $\autlat$ of the top may in turn be identified with $\prod_{1 \le i \le n} \GL\bigl(\bigoplus_{r \in I_i} K \e(r)\bigr)$, where $I_i = \{ r \le t \mid e_i\, \e(r) = \e(r)\}$.   Applying the same considerations to $\aut_{\lahat} (\Phat / \Jhat \Phat)$, we thus obtain a natural isomorphism $\psi: \autlat \rightarrow \autlathat$ of algebraic groups.  Retaining the above notation, our construction yields the following.  (Part (c) is essentially known; cf\. \cite{\AuReSm, Chap. X, Section 2}.) 

\proclaim{Proposition 5.3}  Suppose that $T = \bigoplus_{1 \le i \le n} S_i^{t_i} \in \lamod$ is realizable with respect to the dimension vector $\bd = (d_1, \dots, d_n)$.  Then $\That \in \lahatmod$ is realizable with respect to the dimension vector $\bdhat =   (t_1, \dots, t_n, d_1 - t_1, \dots, d_n - t_n)$ in $(\NN_0)^{2n}$.  Moreover:
\medskip

{\rm \bf{(a)}} There are natural isomorphisms of algebraic groups
$$\autlap / (\autlap)_u \cong \autlat \ \, {\overset \psi \to \cong}\  \autlathat \cong \autlaphat/ (\autlaphat)_u,$$
where $\psi$ is the isomorphism introduced above.  The actions of the unipotent radicals, $(\autlap)_u$ and $(\autlaphat)_u$, on $\grasstbd$ and  $\grass^{\That}_{\bdhat}$, respectively, are trivial, i.e., $f.C = C$ for all $f \in \unirad$, and analogously for the hatted entities.  In particular, the $\autlap$-action on $\grasstbd$ {\rm {(}}resp.,  the $\autlaphat$-action on $\grass^{\That}_{\bdhat}${\rm {)}} is, in fact, an $\autlat$-action {\rm {(}}resp., an $\autlathat$-action{\rm {)}}.   
\medskip

{\rm \bf{(b)}}  The map
$$\Phi^T_{\bd}: \ \grasstbd \rightarrow  \grass^{\That}_{\bdhat}, \ \ \ C \mapsto \Chat$$
is an isomorphism of projective varieties which is equivariant under the group actions in the following sense:  
$$\Phi^T_{\bd}(g.C) = \ghat.\Chat \tag\dagger$$
for $g \in \autlat$, where $\ghat \in \autlathat$ is the image of $g$ under the isomorphism $\psi$.
\medskip

{\rm \bf{(c)}}  The map $\F^T_\bd$ from the set of isomorphism classes of $\la$-modules with dimension vector $\bd$ and top $T$ to the set of isomorphism classes of $\lahat$-modules with dimension vector $\bdhat$ and top $\That$, given by 
$$M = P/C \ \mapsto \ \Mhat = \Phat/ \Chat,$$
is a bijection.  It preserves and reflects direct sum decompositions.  More strongly:  Suppose $M$ is a direct sum of indecomposable submodules with dimension vectors $\bd^{(j)}$ and tops $T^{(j)}$, and $\, \bdhat^{(j)}$ is the $(2n)$-tuple $(\udim T^{(j)}\, , \, \bd^{(j)} - \udim T^{(j)})$ of nonnegative integers.  Then $\Mhat$ is a direct sum of indecomposables with tops $\That^{(j)}$ and dimension vectors $\bdhat^{(j)}$.   Conversely, if $\Mhat$ is a direct sum of indecomposable submodules with dimension vectors $\bdhat^{(j)} = (d^{(j)}_1, \dots, d^{(j)}_{2n})$, then $M$ is a direct sum of indecomposable submodules with dimension vectors 
$$\bd^{(j)} = \bigl( d^{(j)}_1 + d^{(j)}_{n+1}, \dots, d^{(j)}_n + d^{(j)}_{2n} \bigr)$$
and tops $T^{(j)} = \bigoplus_{1 \le i \le n} S_i^{d_i^{(j)}}$.
\medskip

{\rm \bf{(d)}}  For any point $C \in \grasstbd$, the $\la$-submodule lattice of $P/C$ is isomorphic to the $\lahat$-submodule lattice of $\Phat/\Chat$.  If $M'$ is a submodule of $P/C$ with top $T'$ and dimension vector $\bd'$, then the corresponding submodule of $\Phat/\Chat$ belongs to the isomorphism class $\F^{T'}_{\bd'}(M')$.
\endproclaim 

\demo{Proof} (a) The action of $\unirad$ on $\grasstd$ is trivial, because the points of $\grasstd$ are $\la$-submodules of $JP$ and $J^2P = 0$.  In light of the remarks preceding the proposition, the remaining claims under (a) and those under  (b) are straightforward.  Well-definedness of the map $\F^T_{\bd}$ under (c) is a consequence of $(\dagger)$ under (b), and verifying the remaining statements is once more a matter of routine.   \qed \enddemo

\noindent{\bf Remark 5.4.}  The correspondences $\F^T_{\bd}$ of Proposition 5.3 can be pieced together so as to yield a bijection $\F$ from the isomorphism classes of (finitely generated) $\la$-modules to the isomorphism classes of (finitely generated) induced $\lahat$-modules.  Since the only indecomposable $\lahat$-modules which fail to be induced are the simples $\Shat(\ehat_i)$ (see Lemma 5.2), we re-encounter the well-known fact that finiteness of the representation type of $\la$ is equivalent to finiteness of the representation type of $\lahat$ (cf\. \cite{\AuReSm, Chap\. X, Theorem 2.6 }, for instance).  The map $\F$ does not extend to an equivalence between $\lamod$ and the full subcategory of $\lahatmod$ consisting of the induced modules, however, as $\lamod$ and $\lahatmod$ are no more than {\it stably\/} equivalent in general; indeed, compare the endomorphism rings of paired objects $M$ and $\Mhat$ in the presence of loops in $Q$.  

We note  moreover that the image of the restriction of $\F$ to the $\la$-modules with dimension vector $\bd$ contains $\lahat$-modules of differing dimension vectors, depending on the way $\bd$ is split up into dimension vectors of top and radical.  Namely, this image is the union of the isomorphism classes of $\lahat$-modules with dimension vectors $(\udim T, \,\bd - \udim T)$, where $T$ traces the semisimple $\la$-modules which are realizable with respect to $\bd$.  By contrast, if we focus on the restriction of $\F^{-1}$ to the induced $\lahat$-modules of a fixed dimension vector $\bdhat = (d_1, \dots, d_{2n})$, we find the image of this restriction to consist entirely of $\la$-modules with dimension vector 
$$\bold{D} (\bdhat): = \bigl(d_1 + d_{n+1}, \dots, d_n + d_{n + 2n}\bigr).$$
In fact, what makes the passage from $\la$ to $\lahat$ so helpful is the fact that it spreads out the geometric information stored in a single variety $\modlad$ over multiple (irreducible) parameter spaces $\Mod_{\bdhat}(\lahat)$.  
\bigskip

\noindent {\bf{(B) Crossing the bridge:  Generic module properties over algebras with $J^2 = 0$ }} 
\smallskip

Our main focus will be on Kac decompositions of dimension vectors, in particular, on generic indecomposability, and on generic submodule lattices.  The following theorem by Crawley-Boevey and Schr\"oer generalizes results obtained by Kac for hereditary algebras (see \cite{\KacI}, \cite{\KacII}).     

\proclaim{Theorem 5.5} \cite{\CBS, Theorem 1.1}  Let  $\Delta$ be any finite dimensional algebra, $\bd$ a dimension vector for $\Delta$, and $\C$ an irreducible component of $\Mod_\bd(\Delta)$.  Then there exists a {\rm{(}}unique{\rm{)}} decomposition \,$\bd = \sum_{1 \le j \le m} \bd^{(j)}$ of $\bd$ into dimension vectors $\bd^{(j)}$, together with a dense open subset $\U$ of $\C$, such that every module $M$ in $\U$ decomposes in the form 
$M = \bigoplus_{1 \le j \le m} M_j$, where each $M_j$ is indecomposable of dimension vector $\bd^{(j)}$.
\endproclaim

We refer to the above decomposition of $\bd$ as the {\it  Kac decomposition relative to $\C$\/} (suppressing reference to $\Delta$ when there is no danger of ambiguity).   In the situation where $\Delta = KQ$ is hereditary, Schofield provided an algorithm for finding the Kac decomposition of any dimension vector (see \cite{\Scho}).  As we will deduce from Proposition 5.3, this algorithm carries over from $\lahat$ to $\la$.
 
We will say that a dimension vector $\bd'$ is {\it attained on the submodule lattice of a module\/} $M$ if there is a submodule $M' \subseteq M$ with $\udim M' = \bd'$.  Moreover, we will refer to the submodule lattices of the modules in a subvariety $\U \subseteq \modlad$ as the {\it submodule lattices parametrized by\/} $\U$.  By \cite{\BHT, Theorem 4.3},  the (full) sets of dimension vectors attained on the submodule lattices parametrized by $\toptd$ are generically constant.

\proclaim{Theorem 5.6}  Let $\C$ be an irreducible component of $\modlad$, say $\C = \overline{\toptd}$ with $T = \bigoplus_{1 \le i \le n} S_i^{t_i}$, and $\, \bdhat  =   (t_1, \dots, t_n, d_1 - t_1, \dots, d_n - t_n)$.  Then:
\medskip

\noindent {\rm \bf{(a)}}  $\That$ is the generic top of the modules in  $\Mod_{\bdhat}(\lahat)$.  
\medskip

\noindent {\rm \bf{(b)}} The Kac decomposition of $\bd$ relative to $\C$ is determined by the Kac decomposition of $\, \bdhat$ relative to the irreducible variety $\Mod_{\bdhat}(\lahat)$, and vice versa.   More precisely, the former means:  If
$$\bdhat = \sum_{1 \le j \le m} \bdhat^{(j)},$$
is the Kac decomposition of  \, $\bdhat$, and $\bd^{(j)}  = \bold{D}(\bdhat^{(j)})$ in the notation of {\rm \,Remark 5.4}, then $\, \bd = \sum_{1 \le j \le m} \bd^{(j)}$ is the Kac decomposition of $\,\bd$ relative to $\C$.
\smallskip

In particular, the following are equivalent: 

$\bullet$ The modules in $\C$ are generically indecomposable.

$\bullet$ $\bdhat$ is a Schur root of $\Qhat$ {\rm {(}}i.e., $\Mod_{\bdhat}(\lahat)$ contains a module whose endomorphism ring

 \ \ is isomorphic to $K${\rm {)}}.
 
 $\bullet$ Generically, the modules $M$ in $\C$ satisfy $\End_\la(M)/\Hom_\la(M,JM) \cong K$.
 \smallskip

Moreover, the Kac decomposition of $\bdhat$ determines the dimension vectors of the generic tops of the indecomposable summands of the modules in $\C$:  These are precisely the vectors picking out the first $n$ components of the $\bdhat^{(j)}$. 

\medskip

\noindent {\rm \bf{(c)}}  Suppose $\bd = \sum_{1 \le j \le m} \bd^{(j)}$ is the Kac decomposition relative to $\C$, and the $\bd^{(j)}$ are as in {\rm \bf{(b)}}. Then the following conditions are equivalent: 
\smallskip

$\bullet$  $\C$ contains a dense $\GL(\bd)$-orbit.

$\bullet$  $\sum_{1 \le j \le m} \bigl(1 - \langle  \bdhat^{(j)}, \bdhat^{(j)} \rangle\bigr) = 0$, where $\langle - , -\rangle$ is the Euler form of $\Qhat$.
\medskip

\noindent {\rm{\bf(d)}}  The dimension vectors generically attained on the $\lahat$-submodule lattices parametrized by $\Mod_{\bdhat}(\lahat)$ and those generically attained on the $\la$-submodule lattices parametrized by $\C$ are in one-to-one correspondence as follows:  A $(2n)$-tuple $\buhat = (u_1, \dots, u_{2n})$ of nonnegative integers is generically the dimension vector of a submodule of a module in $\Mod_{\bdhat}(\lahat)$ if and only if the $n$-tuple $\bu = (u_1 + u_{n+1}, \dots, u_n + u_{2n})$ is generically the dimension vector of a submodule of a module in $\C$.
\endproclaim

\demo{Proof}  (a)  Since $T$ is realizable with respect to $\bd$, $\That = \allowmathbreak \bigoplus_{1 \le i \le n} \Shat(e_i)^{t_i}$ is realizable with respect to $\bdhat$.  Clearly no simple module of the form $\Shat(e_i)$ occurs in the radical of a $\lahat$-projective cover of $\That$, and therefore any module $\Mhat$ with top $\That$ and dimension vector $\bdhat$ trivially satisfies condition (b) of Proposition 3.9.   Thus, $\That$ is indeed a generic top of $\Mod_{\bdhat}(\lahat)$;  given that $\Mod_{\bdhat}(\lahat)$ is irreducible, $\That$ is {\it the\/} generic top.

Apart from the equivalences under (b) and (c), the remaining statements of the theorem are now immediate consequences of Proposition 5.3.  

As for the equivalences under (b):
Due to Kac, $\bdhat$ is a Schur root of $\Qhat$ if and only if the $\lahat$-modules in $\Mod_{\bdhat} (\lahat)$ are generically indecomposable (see \cite{\KacII, Proposition 1}).  By the first assertion under (b) --  already justified  --  generic indecomposability of the modules in $\Mod_{\bdhat} (\lahat)$ is tantamount to generic indecomposability of the modules in $\C$. 
This shows the first two conditions to be equivalent.  
Moreover: Since, generically, the modules in $\Mod_{\bdhat} (\lahat)$ have top $\That$, the condition that $\bdhat$ be a Schur root means that, generically, $\Mod_{\bdhat} (\lahat)$ consists of modules $\Mhat$ with top $\That$ and trivial endomorphism ring.  Hence we will find the last two conditions to be equivalent as well, provided we can show that the assignment $M \mapsto \Mhat$ in Proposition 5.3(c) is paralleled by the following connection between the endomorphism rings of $M$ and $\Mhat$:
$$\End_{\lahat}(\Mhat)\  \cong\  \End_\la(M)\, / \, \Hom_\la(M,JM).$$

First we observe that the ideal $\Hom_{\lahat} (\Mhat, \Jhat \Mhat)$ of  $\End_{\lahat}(\Mhat)$ is zero, since the top and radical of $\Mhat$ have no simple direct summands in common.  Hence, if $\xhat_1, \dots, \xhat_t$ is the full sequence of top elements of $\Mhat = \Phat/\Chat$, given by $\xhat_r = \e(r) + \Chat$ (see the intoduction to 5(A) for our notation), the subspace $\bigoplus_{1 \le r \le t} K \xhat_r$ is invariant under all endomorphisms of $\Mhat$.  Let $x_1, \dots, x_t$ with $x_r = \e(r) + C$ be the corresponding full sequence of top elements of $M = P/C$.  Then any $\la$-endomorphism $f$ of $M$ can be uniquely written as a sum of an endomorphism $f_1$ which leaves $\bigoplus_{1 \le i \le n} K x_r$ invariant and a map $f_2 \in \Hom_\la(M,JM)$; this uses, once again, the fact that $J^2 = 0$.  By the definition of the correspondence $C \mapsto \Chat$, we thus obtain a $K$-algebra homomorphism $\End_\la(M) \rightarrow \End_{\lahat}(\Mhat)$ which sends $f$ to $\widehat{f_1}$; its kernel is the ideal $\Hom_\la(M,JM)$.   Thus part (b) is proved.

As for the equivalence under (c):  If
$\C$ contains a dense $\GL(\bd)$-orbit, then this orbit is necessarily contained in $\Mod^T_{\bd}$.  Hence Proposition 2.1 guarantees a dense $\autlap$-orbit in $\grasstbd$, and Proposition 5.3(b) provides us with a dense $\aut_{\lahat} (\Phat)$-orbit in $\grass^{\That}_{\bdhat}$.  Clearly, all of these implications are reversible.   But in the $\lahat$-scenario, Kac's Proposition 4 in {\cite{\KacII} tells us that existence of a dense orbit is tantamount to the vanishing of  $\sum_{1 \le j \le m} (1 - \langle  \bdhat^{(j)}, \bdhat^{(j)} \rangle)$.  This completes the argument.  \qed 
\enddemo

\head 6.  Illustrations in the local case.  Generic classification  \endhead

We return to the local case in order to illustrate the transport of information --  regarding Kac decompositions and other generic data  --  from the varieties $\Mod_{\bdhat}(\lahat)$ to the irreducible components of $\modlad$ (Illustration 6.1).

If $\la$ has wild representation type, one would like to at least obtain a classification of the representations in some dense open subset of each irreducible component of $\modlad$, resp., of $\biggrasslad$.  Theoretically, this is possible, provided one does not place any demands of concreteness on the open subset and the modalities of the classification.  Namely, as
was shown by Rosenlicht in \cite{\Ros}, any irreducible variety $X$ which carries a morphic action by an algebraic group $G$ contains a $G$-stable dense open subset which admits a geometric quotient modulo $G$.  However, invoking this existence statement relative to the irreducible components $\C$ of $\modlad$ is of limited value, unless one is able to specify an appropriate open subset of $\C$ in representation-theoretic terms and relate the structure of the encoded modules to the points of the geometric quotient.   

We suggest the following loosely phrased guidelines for a concrete approach to the problem (we believe them to have useful applications only in favorable situations, however):  Structurally describe the modules in a representation-theoretically specified dense open subset of the considered irreducible component; in particular, provide a normal form pinning them down up to isomorphism.  Optimally, the parameters appearing in the normal form trace an algebraic variety $\bold{X}$ which indexes a universal family that makes $\bold{X}$ a fine moduli space for the pertinent modules.  

The class of examples discussed in Illustration 6.1 below provides a good venue for implementing these guidelines; see Generic Classification 6.2.  

As in Section 3(C), we take $\la$ to be 
$$KQ/ \langle \text{all paths of length \ } 2 \rangle,$$
where $Q$ is the quiver with a single vertex and $r$ loops, $\alpha_1, \dots, \alpha_r$. In particular, the dimension vector $\bd$ agrees with its absolute value $d$.  Note: Already in the tame (biserial) case $r = 2$, there are irreducible components $\C$ of $\Mod_d(\la)$ containing indecomposable modules, although, generically, the modules in $\C$ decompose. 

Given a dimension $d\ge 2$, we know from Section 3(C) that 
the irreducible components of $\Mod_d(\la)$ are in $1$-$1$ correspondence with the positive integers $t < d$ such that 
$$t/(d-t) \in [1/r\, , \, r].$$ 
We label these components by the corresponding pairs $(t, d-t)$:  Thus $C_{t, d-t}$ denotes the component of 
$\Mod_d(\la)$ whose modules have generic radical layering $(S^t, S^{d-t})$.

The hereditary algebra $\lahat$ we paired with $\la$ is a generalized Kronecker algebra.  Indeed, $\lahat = K \Qhat$, where $\Qhat$ is the quiver with two vertices, $e_1$, $\ehat_1$, and $r$ arrows from $e_1$ to $\ehat_1$.  By Theorem 5.6, the Kac decomposition of $d$ relative to $\C_{t,d-t}$ is available if we know the Kac decomposition of the dimension vector $\bdhat = (t, d-t)$ for the corresponding Kronecker algebra $\lahat$: Namely, if $\bdhat = \bdhat^{(1)} + \dots +  \bdhat^{(m)} $ with $\bdhat^{(j)} = (d^{(j)}_1, d^{(j)}_2)$ is the Kac decomposition relative to the irreducible variety $\Mod_{\bdhat}(\lahat)$, then $d = d^{(1)} + \cdots + d^{(m)}$ with $d^{(j)}  = d^{(j)}_1 + d^{(j)}_2$ is the Kac decomposition of $d$ relative to the component $\C_{t, d-t}$ of $\modlad$.  We glean additional information from the Kac decpomposition of $\bdhat$.  Namely, let $M_1, \dots, M_m$ represent the indecomposables, of dimension $d^{(j)}$ respectively, that generically arise as direct summands of the modules in $\C_{t, d-t}$.  Then  $\dim \top M_j = d^{(j)}_1$.

For $\Char K = 0$, the Kac decompositions of the dimension vectors over the generalized Kronecker algebras were essentially pinned down in \cite{\KacI, Theorem 4}.  Without any assumption on the characteristic, we exemplify the translation of information about Kronecker algebras, focusing on dimensions $d$ which are congruent to $-1$ modulo $r+1$.  This class of local examples exhibits a wide variety of generic phenomena.

\proclaim{Illustration 6.1}  Suppose $r \ge 2$, $t \ge 1$, and $d - t = r t - 1$.  In particular, $S^t$ is a generic top of $\Mod_d(\la)$.
\smallskip

{\rm {\bf (I) \  $t \le r$.}}  The modules in $\C_{t, rt-1}$ are generically indecomposable.    Moreover, the only indecomposable modules that generically arise as proper submodules of the modules in $\C_{t, rt-1}$ are $S$ and the left regular module $\la$, up to isomorphism.  

{\rm {\bf (I.a)\  $t < r$.}}  The component $\C_{t, rt-1}$ contains infinitely many orbits of maximal dimension.  {\rm (For additional information on the corresponding indecomposable modules, see Generic Classification 6.2.)}

{\rm {\bf (I.b) \  $t = r$.}} The component $\C_{r, r^2-1}$ contains a dense orbit, namely 
that of the module
$$N \ \  = \ \  \biggl(\bigoplus_{1 \le i \le r} \la z_i\ \biggr) \ \bigg/ \ \la \biggl( \sum_{i \le r} \alpha_i z_i\biggr),$$ 
where each $\la z_i$ is a copy of the left regular module. {\rm (Recall that $\alpha_1, \dots, \alpha_r$ are the loops of $Q$.)}
\smallskip

{\rm {\bf (II) \  $t > r$.}} In this case, the modules $M$ in $\C_{t,rt-1}$ decompose generically, in the form
$$M \cong N \oplus \la^{t - r},$$
where $N$ is the generic module for the component $\C_{r, r^2 - 1}$.  In particular, the orbit of $M$ is dense in $\C_{t, rt - 1}$, and the Kac decomposition of $d$ relative to $\C_{t,rt-1 }$ is $d = d^{(0)} + \dots + d^{(t-r)}$ with $d^{(0)} = r + r^2 - 1$ and $d^{(j)} = 1 + r$ for $1 \le j \le t-r$. 
\endproclaim

\demo{Proof}  We consider the dimension vector $\bdhat = (t, rt  - 1)$ for the generalized Kronecker algebra $\lahat$, and denote by $\langle -,- \rangle$ the Euler form of $\Qhat$.   

First let $t \le r$.  To prove the claim concerning generic submodule dimensions of the modules in $\C = \C_{t,rt-1}$, we use Schofield's Theorem 3.2 in {\cite{\Scho}.  It tells us that, generically the $\lahat$-modules with dimension vector $\bdhat$ have a submodule with dimension vector $\bdhat' \le \bdhat$ if and only if $\ext(\bdhat'\, ,\,  \bdhat - \bdhat') = 0$, where $\ext(\ba, \bb) = \min\{ \dim \Ext^1_\la(A,B) \mid A, B \in \lahat\text{-mod},\,  \udim A = \ba,\,  \udim B = \bb\}$.   

Clearly, the only proper indecomposable submodules that occur generically in the $\lahat$-modules with dimension vector $\bdhat$ are either simple, or have dimension vector $(1,r)$, or else a dimension vector $(a, ra - 1)$ for some positive integer $a$ with $a < t$.  To exclude the occurrence of the latter dimension vectors, we compute $\langle (a, ra - 1)\, ,\, \bdhat - (a, ra-1) \rangle$ to be strictly negative, whence, by \cite{\Scho, Theorem 5.4}, we find that $\ext\bigl( (a, ra-1)\, , \, \bdhat -  (a, ra-1) \bigr) > 0$.  Given that generic decomposability of the modules in $\Mod_{\bdhat}(\lahat)$ would amount to a generic summand with dimension vector $(a, ra - 1)$ for some positive $a < t$, we conclude that the $\lahat$-modules with dimension vector $\bdhat$ are generically indecomposable.  Moreover, we glean that the only non-simple indecomposable $\lahat$-module generically arising in the submodule lattices parametrized by $\Mod_{\bdhat}(\lahat)$ is $\lahat$.  Finally, we apply Theorem 5.6 to deduce the corresponding statements for $\la$. 

(I.a). $t < r$. In this case, one finds $1 - \langle \bdhat, \bdhat \rangle$ to be strictly positive.  Consequently, $\C_{t, rt-1}$ fails to contain a dense orbit by Theorem 5.6(c). 
Given that the orbit dimension is lower semi-continuous, irreducibility of $\C_{t,rt-1}$ therefore yields an infinite number of distinct orbits of maximal dimension.  

(I.b). $t=r$.  We compute $\langle (r, r^2-1)\, , \, (r,r^2-1)\rangle = 1$ and apply Theorem 5.6(c) to conclude that $\C_{r,r^2-1}$ contains a dense orbit.  It is, moreover, readily checked that the displayed module $N$ satisfies $\End_\la(N)/\Hom_\la(N,JN) \cong K$ and hence represents the indecomposable module with dense orbit. 

(II).  Now suppose $t >  r$.  As we saw above $(r, r^2 - 1)$ is a Schur root of $\Qhat$, and evidently so is $(1,r)$.  Therefore, verifying the Kac decomposition $\bdhat = \bigl(t\, ,\,  rt - 1\bigr) = \bigl(r\, , \, r^2 - 1\bigr) + (t-r) \cdot \bigl(1\, ,\, r \bigr)$ amounts to showing that $\ext\bigl((r, r^2 - 1)\,, \,  (1,r) \bigr) = 0$\  (the equality $\ext\bigl( (1,r) \, , \, (r, r^2 - 1) \bigr) = 0$ being obvious since the modules with dimension vector $(1,r)$ are generically projective); then \cite{\KacII, Section 4} yields the postulated decomposition of $\bdhat$.  But the vanishing of the relevant generic $\Ext$-dimension again follows from \cite{\Scho, Theorem 5.4}.  These findings for $\lahat$ translate into the claims for $\la$. \qed 
\enddemo

We use the situation addressed in Illustration 6.1 to exemplify the classification goal outlined at the beginning of the section. 
 
\definition{Generic Classification 6.2}  Let $r \ge 2$, $t \ge 1$, and $d = (r + 1)t - 1$.  The irreducible components of $\Mod_d(\la)$ under (I.b) and (II) of Illustration 6.1 contain dense orbits, representing $\la$-modules that we already specified.  

So we assume $t < r$, in which case $\C = \C_{t,rt-1}$ contains infinitely many orbits of maximal dimension.  Generically, the modules in $\C$ have top $T = S^t$ and are of the form  
$$P/ U(\bc) \ \ \ \text{with} \ \ P = \bigoplus_{1 \le i \le t} \la z_i, \ \ \text{and}\ \ U(\bc) = \la \biggl( \sum_{1 \le i \le t} \, \sum_{1 \le j \le r} c_{ij} \alpha_j z_i\biggr),$$
where $\bc = (c_{ij})$ traces the set $\bC$ of those $t \times r$ matrices whose leftmost $t\times t$ minors, $\det \bigl(c_{ij}\bigr)_{1 \le i,j \le t}$, are nonzero.    The isomorphism classes of the listed modules are in bijective correspondence with the following normal forms of the corresponding presentation matrices $\bc \in \bC$:
$$\pmatrix 1& 0 & 0 &\hdots & 0 & c_{1, t+1} &\hdots\  & c_{1,r} \\ 0 & 1 & 0 &\hdots & 0& c_{2, t+1} &\hdots  & c_{2,r} \\  \hdotsfor 8\\  \hdotsfor 8 \\ 0 & 0 & 0 &\hdots & 1 & c_{t, t+1} &\hdots & c_{t,r}  \endpmatrix$$
In fact, if $\bold{X}$ is the subvariety of $\bC$ consisting of the matrices in normal form, then the (informally presented) family $\bigl( P / U(\bc) \bigr)_{\bc \in \bold{X}}$ has the universal property showing $\bold{X}$ to be a fine moduli space for the modules isomorphic to some $P / U(\bc)$ with $\bc \in \bC$.  Clearly, $\dim \bold{X} = t(r-t)$, whence the members of the universal family depend on a non-redundant collection of $t(r-t)$ parameters. 
  \enddefinition 

\demo{Proof of the claims for $t < r$}  Let $T = S^t$ as before, and let $\D$ be the irreducible component of $\grass^T_d$ which corresponds to the irreducible component $\C \cap \Mod^T_d$ of $\Mod^T_d$ under the bijection of Proposition 2.1(b). We observe that the set $\D' = \{U(\bc) \mid \bc \in \bC\}$ is an $\autlap$-stable dense open subset of $\D$.  Basic linear algebra shows that the $\autlap$-orbit (= $\autlat$-orbit) of any point $U(\bc) \in \D'$ consists of the submodules $U(\bw\cdot\bc) \subseteq P$, with $\bw \in \GL_t(K)$.  In fact, on identifying $\autlat$ with $\GL_t$, we readily obtain an equivariant isomorphism of varieties $\D' \cong \bC$.  In particular, this shows every $\autlat$-orbit of $\D'$ to contain precisely one element in normal form.  Moreover, we see that a geometric quotient of $\D'$ by $\autlat$, if existent, coincides with a geometric quotient of the variety $\bC$ by its left $\GL_t$-action.  To confirm existence, one checks that the morphism $\bC \rightarrow \bold{X}$, which sends any matrix $( \ba | \bb)$ in $\bC$ (where $\ba$ has size $t\times t$) to the matrix $\ba^{-1} (\ba | \bb)$ in $\bold{X}$, is indeed a geometric quotient of $\bC$ by $\GL_t$.  It is now routine to verify that the obvious bundle of $\la$-modules in $\D'$ which is parametrized by $\bold{X}$  --  it formalizes the family $\bigl( P / U(\bc) \bigr)_{\bc \in \bold{X}}$  -- satisfies the universal property making $\bold{X}$ a fine moduli space for the modules in $\D'$ (cf\. \cite{\moduli}). \qed 
\enddemo

\head 7.  The dense orbit property for $J^2 = 0$ \endhead

In this section, we give another application of the geometric link, exhibited in Section 5, between algebras with vanishing radical square and hereditary algebras. 

It is readily verified that any finite dimensional algebra $\la$ 
of finite representation type satisfies the following {\it dense orbit property} (terminology of \cite{\CKW}):  Namely, for any dimension vector $\bd$, each of the irreducible components of $\Mod_{\bd}(\la)$ is the closure of a single orbit.  Chindris, Kinser and Weyman posed the following problem:  For which classes of algebras does the converse hold as well, i.e., for which algebras does the dense orbit property imply finite representation type? While this is well known to be true for hereditary algebras, they demonstrated failure in general.  Among the classes of algebras for which they showed the converse to be true  are the string algebras and the algebras admitting a preprojective component (see \cite{\CKW, Sections 3, 4}).

We will invest the geometric $Q$-$\Qhat$ connection of Section 5 to add the algebras with vanishing radical square to the list of positive instances.  For our notation, we refer to Section 5.

\proclaim{Lemma 7.1}  Let $J^2 = 0$, and suppose that $\bdhat = (d_1, \dots, d_{2n})$ is a dimension vector such that $\Mod_{\bdhat}(\lahat)$ contains an indecomposable induced module $\Mhat$.  If $\Mhat$ has top $\That = \bigoplus_{1 \le i \le n} (\Shat(e_i))^{t_i}$ and $\bd = (d_1+d_{n+1}, \dots, d_n + d_{2n})$, then $T =  \bigoplus_{1 \le i \le n}S_i^{t_i}$ is a generic top of $\modlad$.
\endproclaim

\demo{Proof}  Suppose that $\Mhat$ is as in the claim, say $\Mhat \cong \Phat/\Chat$ with $\Chat \in \grass^{\That}_{\bdhat}$.  We apply Proposition 5.3 to deduce that then $M = P/C$ (where $P$ is the projective cover of $T$ and $C = \bigl(\Phi^{T}_{\bd} \bigr)^{-1}(\Chat) \in \grasstbd$) is an indecomposable module in $\toptd$.  Either $M$ is simple, in which case $T=M$ is trivially a generic top of $\modlad$, or else $M$ is without simple direct summands, in which case Corollary 3.7 implies that $T$ is a generic top of $\modlad$.  \qed \enddemo

\proclaim{Theorem 7.2} Suppose $J^2 = 0$.  Then the following conditions are equivalent:
\smallskip

\noindent{\rm \bf{(1)}}  $\la$ has the dense orbit property.
\smallskip

\noindent{\rm \bf{(2)}} For every dimension vector $\bd$, any irreducible component of $\modlad$ which contains an indecomposable module has a dense orbit.
\smallskip

\noindent{\rm \bf{(3)}}  $\la$ has finite representation type.
\endproclaim

\demo{Proof}  The implication ``(3)$\implies$(1)" is known, and ``(1)$\implies$(2)" is trivial.  To prove that (2) implies (3), suppose that (2) holds.  We will deduce that then also $\lahat$ satisfies (2), meaning that every module variety $\Mod_{\bdhat}(\lahat)$ which contains an indecomposable $\lahat$-module $N$ has a dense orbit.  Since this is clear when $N$ is simple, we assume that $N$ in $\Mod_{\bdhat}(\lahat)$ is indecomposable, but not simple.  Thus, by Lemma 5.2, $N$ is induced say $N = \Mhat$.  Let $\That = \top \Mhat$, and use Lemma 7.1 to ascertain that the corresponding semisimple $\la$-module $T$ is a generic top of $\modlad$, where $\bd = (d_1 + d_{n+1}, \dots, d_n + d_{2n})$ is the matching dimension vector of $\la$.  Let $M$ in $\toptd$ represent the isomorphism class assigned to that of $\Mhat$ by the bijection $\F^T_{\bd}$ of Proposition 5.3(c).  Then this proposition shows $M$ to be indecomposable as well, whence our hypothesis guarantees a dense orbit in $\C = \overline{\toptd}$.  Now Theorem 5.3(b) yields a dense orbit in $\grass^{\That}_{\bdhat}$. 

It is well-known that the analogue of condition (2) for $\lahat$ forces $\lahat$ to have finite representation type.  Since we were unable to locate a reference providing exactly what we need, we include the short argument. Again, we denote by $\langle -\, , \,  - \rangle$ the Euler form of $\lahat$.  Let $\bdhat$ be any dimension vector of $\lahat$, and let $\bdhat = \sum_{1 \le j \le m} \bdhat^{(j)}$ be its Kac decomposition.  This means that, generically, the modules $N$ in $\Mod_{\bdhat}(\lahat)$ decompose in the form $N = \bigoplus_{1 \le j \le m} N_j$, where the $N_j$ have the following properties:  $N_j$ is indecomposable with dimension vector $\bdhat^{(j)}$, such that the $\GL(\bdhat^{(j)})$-orbit, $\orbit(N_j)$, that corresponds to the isomorphism class of $N_j$ in $\Mod_{\bdhat^{(j)}}(\lahat)$ has maximal dimension, and $\Ext^1_{\lahat}(N_i , N_j) = 0$ for $i \ne j$.
Since $\lahat$ satisfies (2), each of the varieties $\Mod_{ \bdhat^{(j)} } (\lahat)$ contains a dense orbit, and consequently, the direct summands $N_j$ in the above generic decomposition are unique up to isomorphism.  Therefore $\dim \orb(N_j) = \dim \Mod_{ \bdhat^{(j)} } (\lahat)$, which in turn implies $\langle \bd^{(j)}, \bd^{(j)}\rangle$ to be positive.  Due to the vanishing of the mixed $\Ext$-spaces, we conclude
$$\langle \bdhat, \bdhat \rangle = \dim \bigl(\bigoplus_{j=1}^m  N_j \bigr)\ - \ \dim \Ext^1_{\lahat}\bigl(\bigoplus_{j=1}^m  N_j , \bigoplus_{j=1}^m  N_j \bigr)\ \  =\ \ \sum_{j=1}^m \langle \bdhat^{(j)},  \bdhat^{(j)}\rangle \ \ > \ \ 0.$$
Thus the Tits form of $\lahat$ is positive definite, meaning that $\lahat$ has finite representation type.

Finally, we use the well-known fact that finite representation type is passed on from $\lahat$ to $\la$ (see Remark 5.4) to find that $\la$ indeed satisfies (3). \qed
\enddemo

\enddocument